\documentclass[12pt,onecolumn,twoside]{IEEEtran}
\usepackage{amsmath,amssymb,wasysym,yfonts,psfrag,latexsym,dsfont,graphicx,bbm} 

\begin{document}
\def\spacingset#1{\def\baselinestretch{#1}\small\normalsize}

\newcounter{plm}
\newcounter{plm1}
\newcounter{plm2}
\newcounter{plm3}
\newcounter{plm4}
\newcounter{plm5}
\newcounter{plm6}
\newcounter{plm7}
\newtheorem{thm}[plm]{Theorem}
\newtheorem{lemma}[plm1]{Lemma}
\newtheorem{cor}[plm2]{Corollary}
\newtheorem{remark}[plm3]{Remark}
\newtheorem{prop}[plm4]{Proposition}
\newtheorem{problem}[plm5]{Problem}
\newtheorem{ex}[plm6]{Example}
\newtheorem{definition}[plm7]{Definition}

\newcommand{\mR}{{\mathbb R}}
\newcommand{\mZ}{{\mathbb Z}}
\newcommand{\mN}{{\mathbb N}}
\newcommand{\mM}{{\mathbb M}}
\newcommand{\mH}{{\mathbb H}}
\newcommand{\mS}{{\mathbb S}}
\newcommand{\mI}{{\mathbb I}}
\newcommand{\mE}{{\mathbb E}}
\newcommand{\mC}{{\mathbb C}}
\newcommand{\mU}{{\mathbb U}}
\newcommand{\mmu}{{\mathbb u}}

\newcommand{\fR}{{\mathfrak R}}
\newcommand{\fG}{{\mathfrak G}}
\newcommand{\fT}{{\mathfrak T}}
\newcommand{\fC}{{\mathfrak C}}
\newcommand{\fF}{{\mathfrak F}}
\newcommand{\fM}{{\mathfrak M}}

\newcommand{\cM}{{\mathcal M}}
\newcommand{\cS}{{\mathcal S}}
\newcommand{\cR}{{\mathcal R}}
\newcommand{\cU}{{\mathcal U}}
\newcommand{\cG}{{\mathcal G}}
\newcommand{\cC}{{\mathcal C}}
 \newcommand{\cN}{{\mathcal N}}
 \newcommand{\cK}{{\mathcal K}}
 \newcommand{\cI}{{\mathcal I}}
 
 \newcommand{\bR}{{\bf R}}
 
 \newcommand{\rleft}{{\rm left}}
 \newcommand{\rright}{{\rm right}}
 \newcommand{\trace}{{\rm trace}}
 \newcommand{\symmetric}{{\rm Herm}}
 \newcommand{\Herm}{{\rm Herm}\,}
\newcommand{\Real}{{\Re}e\,}
\newcommand{\la}{{z}}

\newcommand{\eGamma}{{\it \Gamma}}
\newcommand{\cL}{{\cal L}}
\newcommand{\Cp}{\eGamma}
\newcommand{\mD}{{\mathbb D}}
\newcommand{\cH}{{\cal H}}
\newcommand{\cW}{{\cal W}}
\newcommand{\cV}{{\mathcal V}}
\newcommand{\cF}{{\mathbb{F}}}
\newcommand{\cE}{{\mathcal E}}
\newcommand{\bPi}{{\bf \Pi}}

\title{The Carath\'{e}odory-Fej\'{e}r-Pisarenko decomposition
and its multivariable counterpart}
\author{Tryphon T. Georgiou, {\em Fellow, IEEE} \thanks{Department of Electrical and Computer Engineering, University of Minnesota, Minneapolis, MN 55455; {\tt tryphon@ece.umn.edu}
\hspace*{10pt}Research partially supported by the NSF and the AFOSR.}}
\date{}
\markboth{Submitted to IEEE Trans.\ on Automatic Control
}
{Georgiou: The Carath\'{e}odory-Fej\'{e}r-Pisarenko decomposition}
\maketitle

\begin{abstract} When a covariance matrix with a Toeplitz structure is written as the sum of a singular one and a positive scalar multiple of the identity, the singular summand corresponds to the covariance of a purely deterministic component of a time-series whereas the identity corresponds to white noise---this is the Carath\'{e}odory-Fej\'{e}r-Pisarenko (CFP) decomposition. In the present paper we study multivariable analogs for block-Toeplitz matrices as well as for matrices with the structure of state-covariances of finite-dimensional linear systems (which include block-Toeplitz ones). We characterize state-covariances which admit only a deterministic input power spectrum. We show that multivariable decomposition of a state-covariance in accordance with a ``deterministic component $+$ white noise'' hypothesis for the input does not exist in general, and develop formulae for spectra corresponding to singular covariances via decomposing the contribution of the singular part. 
We consider replacing the ``scalar multiple of the identity'' in the CFP decomposition by a covariance of maximal trace which is admissible as a summand. The summand can be either (block-)diagonal corresponding to white noise or have a ``short-range correlation structure'' correponding to a moving average component. 
The trace represents the maximal variance/energy that can be accounted for by a process (e.g., noise) with the aforementioned structure at the input, and the optimal solution can be computed via convex optimization. The decomposition of covariances and spectra according to the range of their time-domain correlations is an alternative to the CFP-dictum with potentially great practical significance.
\end{abstract}


\section{\bf Introduction}
\PARstart{P}{resent} day signal processing is firmly rooted in the analysis and interpretation of second order statistics. In particular, the observation that
singularities in covariance matrices reveal a deterministic linear dependence between observed quantities, forms the basis of a wide range of techniques, from Gauss' least squares to modern subspace methods in time-series analysis. In the present work we study the nature and origin of singularities in certain structured covariance matrices which arise in multivariable time-series.

Historically, modern subspace methods (e.g., MUSIC, ESPRIT) can be traced to Pisarenko's harmonic decomposition and even earlier to a theorem by C.\ Carath\'{e}odory and L.\ Fej\'{e}r on a canonical decomposition of finite Toeplitz matrices \cite{GrenanderSzego,Haykin,StoicaMoses}. The Toeplitz structure characterizes covariances of stationary scalar time-series. Their multivariable counterpart, block-Toeplitz matrices, having a less stringent structure, has received considerably less attention. The present work focuses on analogues of the Carath\'{e}odory-Fej\'{e}r-Pisarenko (CFP) decomposition to finite block-Toeplitz matrices as well as to the more general setting of state-covariances of a known linear dynamical system.

In Section \ref{sec1} we begin with background material on matrices with the structure of a state-covariance of a known linear dynamical system---block-Toeplitz matrices being a special case. Section \ref{analyticinterpolation} discusses the connection between covariance realization and analytic interpolation. Section \ref{dualformalism} presents a duality between left and right matricial Carath\'{e}odory interpolation and their relation to the time arrow in dynamical systems generating the state-process. Duality is taken up again in Section \ref{predictionpostdiction} where we study optimal prediction and postdiction (i.e., prediction backwards in time) of a stochastic input based on state-covariance statistics. The variance of optimal prediction and postdiction errors coincide with left and right uncertainty radii in a Schur representation of the family of consistent spectra given in \cite{acmatrix1,acmatrix2} and elucidate the symmetry observed in these references. Further, Section \ref{predictionpostdiction} presents geometric conditions on the state-covariance for the input process to be deterministic and for the optimal predictor and postdictor to be uniquely defined. Vanishing of the variance of the optimal prediction or postdiction errors is shown in Section \ref{sec:singleton} to characterize state-covariances for which the family of consistent input spectra is a singleton.

Section \ref{positiverealfunction} gives a closed form expression for the power spectrum
corresponding the ``central solution'' of \cite{acmatrix2}. This result extends the theory in  \cite{acmatrix2} to the case where the state-covariance is singular. Naturally, the subject of this section has strong connections with the theory of Szeg\"{o}-Geronimus orthogonal polynomials and their multivariable counterparts \cite{dgk1}. In this section, we present yet another generalization of such polynomials as they now become matricial functions sharing the eigen-structure of the transfer function of the underlying dynamical system. Then, Section \ref{residues} explains how to isolate the deterministic component of the power spectrum via computation of relevant residues with matrix techniques.

Section \ref{whitepluscolor} shows, by way of example, that a state-covariance may not admit a decomposition into one corresponding to white-noise plus another corresponding to a deterministic input. To this end, a natural generalization of the
CFP decomposition is to seek a maximal white-noise component at the input consistent with a given state-covariance. We explain how this is computed and discuss yet a further generalization where the input ``noise'' is allowed to have ``short-range correlation structure''. 
For instance, if the state-covariance is $\ell\times \ell$ (block-)Toeplitz, then we may seek to account for input noise whose auto-covariance vanishes after the $k<\ell$-moment---i.e., colored noise modeled by at most a $k$-order moving average filter. In this way, a maximal amount of variance that may be due to short range correlations can be accounted for, leaving the remaining energy/variance to be attributed to periodic deterministic components and possibly, stochastic components with long range (longer than $k$) correlations.

\section{\bf Structured covariance matrices}\label{sec1}

Throughout we consider a multivariable, discrete-time, zero-mean, stochastic process
\[
\{u_k\;:\; k\in\mZ\}
\]
taking values in $\mC^{m\times 1}$ with $m\in\mN$. Thus, $u_k$ is to be thought of as a column vector. We denote by
\[
R_k:=\cE\{u_\ell u_{\ell-k}^*\},
\]
for $k,\ell\in\mZ$, the sequence of matrix covariances and by
$d\mu(\theta)$ the corresponding matricial spectral measure for which
\[
R_k=\int_0^{2\pi} e^{-jk\theta}d\mu(\theta)
\]
for $k\in\mZ$ (see e.g., \cite{masani}). As usual, star ($^*$) denotes the complex-conjugate transpose of, prime ($^\prime$) denotes the transpose, $j:=\sqrt{-1}$ following the usual ``engineering'' convention, and $\cE\{\cdot\}$ denotes the expectation operator. Whenever star ($^*$) is applied to a rational function of $\la$ it represents the para-conjugate Hermitian $f(\la)^*:=f^*(\la^{-1})$ where $f^*(\cdot)$ refers to $*$-ing the coefficients of $f(\cdot)$ whereas the transformation of the argument is indicated separately.

It is well-known that a covariance sequence
\[
\{R_\ell\;:\;\;\ell\in\mZ\mbox{ and }R_{-\ell}=R_\ell^*\}
\]
is completely characterized by the non-negativity of the block-Toeplitz matrices
\begin{equation}\label{bRsubn}
\bR_\ell:=\left[ \begin{array}{cccc} R_0 & R_1 &\ldots & R_\ell\\
                                       R_{-1} & R_0 & \ldots & R_{\ell-1}\\
                                       \vdots & \vdots & \ddots &\vdots\\
                                       R_{-\ell} & R_{-\ell+1} & \ldots &R_0\end{array}\right]
\end{equation}
for all $\ell$. That is, such an infinite sequence with the property that $\bR_\ell\geq 0$, $\forall \ell$, qualifies as a covariance sequence of a stochastic process and vice versa. On the other hand, the infinite sequence of $R_\ell$'s defines the spectral measure $d\mu$ (up to an additive constant) and conversely.

It is often the case that only a finite set of second-order statistics is available, and then, it is of interest to characterize possible extensions of the finite covariance sequence $\{R_0,R_1,\ldots,R_\ell\}$, or equivalently, the totality of consistent spectral measures (see \cite{dgk1,dgk2,dgk3,DD,BGL2,acmatrix1,acmatrix2}). In general, these are no longer specified uniquely by the finite sequence $\{R_0,R_1,\ldots,R_\ell\}$.
In the present paper we are interested in particular, in the case where a finite set of second-order statistics such as $\{R_0,R_1,\ldots,R_\ell\}$ completely specifies the corresponding spectral measure (and hence, any possible infinite extension as well). We address this question in the more general setting of structured covariance matrices which includes block-Toeplitz matrices as a special case.

A block-Toeplitz matrix such as $\bR_\ell$ given in (\ref{bRsubn}) can be thought of as the state-covariance of the linear (discrete-time) dynamical system
\begin{eqnarray}
x_k&=&A x_{k-1} + B u_k,                                  \label{i2s}
\mbox{ for $k\in{\mathbb Z}$}.
\end{eqnarray}
where
\begin{equation}\label{companion}
A=\left[ \begin{array}{ccccc} O_m & O_m &\ldots & O_m& O_m\\
                                       I_m & O_m & \ldots & O_m& O_m\\
                                        & \ddots & \ddots &\vdots&\vdots\\ \\
                                       O_m & O_m &  &I_m& O_m\end{array}\right], 
                                       B=\left[ \begin{array}{c}
                                       I_m\\ O_m\\ \vdots\\ \\ O_m\end{array}\right]
\end{equation}
with $O_m$ and $I_m$ the zero and the identity matrices of size $m\times m$, $A$ a $(\ell+1)\times(\ell+1)$ and $B$ a $(\ell+1)\times 1$ block matrices, respectively. The size of each block is $m\times m$ and hence the actual sizes of $A,B$ are $n\times n$ and $n\times m$, with $n=(\ell + 1)m$, respectively.
While for general state-matrices $A,B$ the structure of the state-covariance may not be visually recognizable, it is advantageous, for both, economy of notation and generality, to   develop the theory in such a general setting---the theory of block-Toeplitz matrices being a special case.

Thus, henceforth, we consider an {\em input-to-state} dynamical system as in (\ref{i2s}) where
\begin{eqnarray}\nonumber
&(\ref{standing1}a)&u_k\in\mC^m,\; x_k\in\mC^n,\; A\in\mC^{n\times n},\; B\in\mC^{n\times m},\\\nonumber
&(\ref{standing1}b)&{\rm rank}(B)=m,\\\nonumber
&(\ref{standing1}c)&(A,B) \mbox{ is a reachable pair, and}\\\nonumber
&(\ref{standing1}d)&\mbox{all the eigenvalues of $A$}\\
&&\mbox{have modulus } < 1. \label{standing1}
\end{eqnarray}
Without loss of generality and for convenience we often assume that the pair $(A,B)$ has been normalized as well so that\\[.02in]
\[
\hspace*{-4cm}(\ref{standing1}e)\;\;\,AA^*+BB^*=I_n.
\]
Conditions (\ref{standing1}a-d) are standing assumptions throughout. Whenever condition (\ref{standing1}e) is assumed valid, this will be stated explicitely.
With $u_k\in\mC^m$, $k\in\mZ$, a
zero-mean stationary stochastic process we denote by
\[
\bR:=\cE\{x_kx_k^*\}
\]
the corresponding (stationary) state-covariance. The space of Hermitian $n\times n$ matrices will be denoted by $\mH_n\subset\mC^{n\times n}$ while positive (resp.\ nonegative) definiteness of an $\bR\in\mH_n$ will be denoted by $\bR>0$ (resp.\ $\bR\geq 0$). Any state-covariance as above certainly satisfies both conditions, i.e., it is Hermitian and non-negative definite. The following statement characterizes the linear structure imposed by (\ref{i2s}).\\

\begin{thm}\label{thm1} {\em (see \cite{acmatrix1}):} {\sf A nonnegative-definite Hermitian matrix $\bR$ (i.e., $\mH_n\ni\bR\geq 0$) arises as the (stationary) state-covariance of (\ref{i2s}) for a suitable stationary input process $\{u_k\}$ if and only if the following equivalent conditions hold:
\begin{eqnarray}\nonumber
(\ref{rankcondition}a)&&
{\rm rank}\left[\begin{array}{cc}\bR-A\bR A^*&B\\B^*&0\end{array}\right]
=
2m,\\
&&\mbox{or, equivalently,}\nonumber\\ && \nonumber\\ (\ref{rankcondition}b)&&
\bR-A\bR A^*=BH+H^*B^*\nonumber\\
&&\mbox{for some }H\in\mC^{m\times n}.\label{rankcondition}
\end{eqnarray}
}\end{thm}
\vspace*{.1in}

\begin{proof} See \cite[Theorems 1 \& 2]{acmatrix1}.\end{proof}
\vspace*{.1in}

A finite $m\times m$ non-negative matrix-valued measure $d\mu(\theta)$ with $\theta\in(-\pi,\pi]$
represents the power spectrum of a stationary $m\times 1$-vector-valued stochastic process.
The class of all such $m\times m$ matrix-valued non-negative bounded measures will be denoted by $\mM$. Note that the size $m$ is suppressed in the notation because it will be the same throughout. Starting with a stationary input $u_k$ with power spectral distribution
$d\mu\in \mM$, the state-covariance of (\ref{i2s}) can be expressed in the form of the integral (cf.\ \cite[Ch.\ 6]{masani})
\begin{equation}\label{Sigma}
\bR=\int_0^{2\pi}
\left( G(e^{j\theta}) \frac{d\mu(\theta)}{2\pi} G(e^{j\theta})^* \right)
\end{equation}
where
\[
G(\la):=(I_n-\la A)^{-1}B
\]
is the transfer function of (\ref{i2s}) (with $z$ corresponding to the delay operator, so that ``stability'' corresponds to ``analyticity in the open unit disc $\mD:=\{\la\in\mC : |\la|<1\}$''). Thus, either condition (\ref{rankcondition}a) or (\ref{rankcondition}b) in the above theorem characterizes the range of the mapping
\[
\mM\ni d\mu \mapsto \bR
\]
specified by (\ref{Sigma}). The family of power spectral distributions which satisfy (\ref{Sigma}) will be denoted by
\[
\mM_\bR:=\{ d\mu(\theta)\in\mM\;:\; \mbox{ equation (\ref{Sigma}) holds}\}.
\]
The above theorem states that this family is nonempty when $\bR$ satisfies the stated conditions. Furthermore,
a complete parametrization of $\mM_\bR$ is given in \cite{acmatrix1,acmatrix2}.

The present work explores the case where $\mM_\bR$ is a singleton. The special
case where $u_k$ is scalar and $\bR$ a Toeplitz matrix (but not ``block-Toeplitz'')
goes back to the work of Carath\'{e}odory and Fej\'{e}r a century ago, and later on, to the work of Pisarenko (see \cite{GrenanderSzego,Haykin,StoicaMoses}). In the scalar case, $\mM_\bR$ is a singleton if and only if $\bR$ is singular (and of course non-negative definite). Then  $u_k$ is deterministic with a spectral distribution $d\mu$ having at most $n-1$ discontinuities (spectral lines). In the present paper we obtain analogous results when $\bR$ is a state-covariance and $\mM_\bR$ is a singleton, and then we study decomposition of a general $\bR>0$ into a covariance due to ``noise'' plus a singular covariance with deterministic components---in the spirit of the
CFP 
decomposition of Toeplitz covariance matrices.

\section{\bf Connection with analytic interpolation}\label{analyticinterpolation}

The early work of Carath\'{e}odory and Fej\'{e}r was motivated by questions in analysis which led to the development of analytic interpolation theory---a subject which has since attained an important place in operator theory, and more recently, closer to home, in robust control engineering. We review certain rudimentary facts and establish notation.

A non-negative measure $\mu\in\mM$ 
specifies an $m\times m$ matrix-valued function
\begin{eqnarray}\nonumber
F(\la)&=&
\int_0^{2\pi}\left(\frac{1+\la e^{j\theta}}{1-\la e^{j\theta}}\right)\frac{d\mu(\theta)}{2\pi} + jc,\\
&=:& \cH[d\mu] +jc\label{correspondingpr}
\end{eqnarray}
with $jc$ an arbitrary skew-Hermitian constant (i.e., $c\in\mH_m$),
which is analytic in the open unit disc $\mD$
and has non-negative definite Hermitian part (see, e.g., \cite[page 36]{dgk3}). We denote by $\cH[d\mu]$ the {\em Herglotz} integral given in previous line. The class of such
$m\times m$ functions with non-negative Hermitian part in $\mD$, herein denoted by
\begin{eqnarray*}
\cF&:=&\{F(\la)\;:\;F(\la)=\cH[d\mu]+jc\\
&& \mbox{ with }c\in\mH_m\mbox{ and }\mu\in\mM\},
\end{eqnarray*}
is named after Carath\'{e}odory and often referred to simply as ``positive-real''.
Conversely, given $F\in\cF$, a corresponding $d\mu(\theta)$ can be recovered by the radial (weak) limits of the Hermitian part of $F(\la)$; 
\begin{equation}\label{realpart}
d\mu(\theta)= \lim_{r\nearrow 1}{\Herm}\{F(re^{j\theta})\}.
\end{equation}
In fact these two families, $\cF$ and $\mM$, are in exact correspondence via (\ref{correspondingpr}) and (\ref{realpart}) (assuming that elements in $\cF$ are identified if they only differ by a skew-Hermitian constant and, similarly, non-decreasing distribution functions $\mu$  are defined up to an arbitrary additive constant).

Given $(A,B)$ as above, let
$C\in\mC^{m\times n}$, $D\in\mC^{m\times m}$ be selected so that
\begin{equation}\label{V}
V(\la):=D+\la C(I_n-\la A)^{-1}B
\end{equation}
is inner, i.e., $V(\xi)^*V(\xi)
=I_m$ for all $|\xi|=1$. Since
$V(\la)$ is square, $V(\xi)V(\xi)^*=I_m$ as well. If the normalization (\ref{standing1}e) is in place, the condition on $(C,D)$ for $V(\la)$ to be inner is simply that
\[
U:=\left[\begin{array}{cc}A&B\\C&D\end{array}\right]
\]
is a unitary matrix.
The rows of $G(\la)$ form a basis of
\[
\cK:=\cH_2^{1\times m}\ominus\cH_2^{1\times m}V(\la)
\]
where $\cH_2$ denotes the Hardy space of functions analytic in $\mD$ with square-integrable boundary limits. This can be easily seen from the identity \cite[equation (38)]{acmatrix1}
\begin{equation}\label{represAC}
G(\la)=(\la I-A^*)^{-1}C^*V(\la)
\end{equation}
(from which it follows that the entries of $G(\la)V(\la)^*$ are in $\cH_2^\perp$, the orthogonal complement of $\cH_2$ in the Lebesgue space of square-integrable function on the unit circle $\cL_2(\partial \mD)$).

Now let $d\mu(\theta)$ represent the power spectrum of the input to (\ref{i2s}),
$\bR$ the corresponding state-covariance, and $F(\la)$ obtained via (\ref{correspondingpr}).
Then, $\bR$ turns out to be the Hermitian part of the operator
\begin{equation}\label{cW}
\cW\;:\;\cK\to \cK : \nu(\la)\mapsto \bPi_\cK  \left(\nu(\la) F(\la)^*\right),
\end{equation}
with respect to basis elements being the rows of $G(\la)$,
where $\bPi_\cK$ denotes the orthogonal projection onto $\cK$ (see \cite[equations (40-41)]{acmatrix1}). Of course, $\bR$ is also the Grammian
\[
\langle G(\la),G(\la)\rangle_{d\mu}
\]
with respect to the inner product
\[
\langle g_i(\la),g_k(\la)\rangle_{d\mu}:=\int_0^{2\pi}\left(g_i(e^{j\theta})\frac{d\mu(\theta)}{2\pi}g_k(e^{j\theta})^*\right).
\]
This is in fact the content of (\ref{Sigma}).

The relationship between $F(\la)$ and $\bR$ can be obtained by way of $\cW$.
If $H$ is the zeroth Fourier coefficient of $G(\la)F(\la)^*$ then the matrix representation $W$
for $\cW$ with respect to the rows of $G(\la)$ satisfies (see \cite{acmatrix1})
\begin{equation}\label{Wequation}
W-AWA^*=H^*B^*
\end{equation}
leading to (\ref{rankcondition}) for $\bR=W+W^*$.
The matrices $W$ or $H$ completely specify $\bPi_\cK F(\la)^* |_\cK$ and in fact
\begin{equation}\label{PRNehari}
F(\la)=F_0(\la) + Q(\la) V(\la)
\end{equation}
with
\[
F_0(\la):=H(I_n-\la A)^{-1}B
\]
and $Q(\la)$ is a matrix-valued function which is analytic in $\mD$.
Conversely, if $F(\la)\in\cF$ and satisfies
(\ref{PRNehari}), then it gives rise via (\ref{realpart}) to a measure
which is consistent with the state-covariance $\bR$.

Equation (\ref{PRNehari}) specifies a problem which is akin to the Nehari problem encountered in
${\cal H}_\infty$-control theory, but involves interpolation with
positive-real functions instead of functions in ${\cal H}_\infty(\mD)$.
Some of the early work in analytic interpolation focused on conditions in terms of interpolating
values $F(\la_i)$ at specified points $\la_i\in\mD$ ($i=1,\ldots,n$) which guarantee the existence of a scalar $F(\la)\in\cF$. Invariably, the conditions involve the non-negativity of the so-called Pick matrix. In the current setting the corresponding Pick matrix is non other than $\bR$ (see \cite{acmatrix1,acmatrix2}). For further references and trends in literature on analytic interpolation see \cite{DD,bgr}.

\section{\bf A dual formalism}\label{dualformalism}

Using (\ref{represAC}), equation (\ref{Sigma}) can be rewritten as
\begin{equation}\label{Sigmar}
\bR=\int_0^{2\pi}
\left( G_r(e^{j\theta}) \frac{d\mu_r(\theta)}{2\pi} G_r(e^{j\theta})^* \right),
\end{equation}
where
\[
G_r(\la)=(\la I-A^*)^{-1}C^*
\]
and
\begin{equation}\label{spectralmeasures}
d\mu_r(\theta)=V(e^{j\theta})d\mu(\theta)V(e^{j\theta})^*.
\end{equation}
The rows of $G_r(\la)$, for $\la=e^{j\theta}$, span a subspace of $\left(\cH_2^{1\times m}\right)^\perp$ which we denote by
\[
\cK_r:=\left(\cH_2^{1\times m}\right)^\perp\ominus \left(\cH_2^{1\times m}\right)^\perp V(\la)^*.
\]
The notation $^\perp$ denotes orthogonal complement in the ``ambient'' space---here $\cL_2(\partial \mD)^{1\times m}$.
It readily follows that a state-covariance of (\ref{i2s}) satisfies a set of dual conditions
given below.\\

\begin{thm}\label{thm2} {\sf A nonnegative-definite Hermitian matrix $\bR\in\mC^{n\times n}$ arises as the (stationary) state-covariance of (\ref{i2s}) for a suitable stationary input process $u_k$ if and only if the following equivalent conditions hold:
\begin{eqnarray}\nonumber
(\ref{rankcondition}c)&&
{\rm rank}\left[\begin{array}{cc}\bR-A^*\bR A&C^*\\C^*&0\end{array}\right]
=
2m,\\
&& \mbox{or, equivalently,}\nonumber\\
(\ref{rankcondition}d)&&
\bR-A^*\bR A=C^*L^*+LC\nonumber\\ &&\mbox{for some }L\in\mC^{n\times m}\nonumber
\end{eqnarray}
and $C$ selected as in Section \ref{analyticinterpolation} (i.e., so that $D+zC(I_n-zA)^{-1}B$ is inner). Conditions (\ref{rankcondition}c-d) are also equivalent to conditions (\ref{rankcondition}a-b).
}\end{thm}
\vspace*{.1in}

\newcommand{\rank}{{\rm rank}}
It is noted that, $\rank(B)=m$ in condition (\ref{standing1}b) implies that
$\rank(C)=m$ as well. To see this, assume without loss of generality that (\ref{standing1}e) holds. Then $B^*B=I_m-D^*D>0$ which implies that $\|D\|<1$. Using once more unitarity of $U$ and the fact that $\|D\|<1$, we obtain that $CC^*=I_m-DD^*>0$ which implies that $\rank(C)=m$.

An insightful derivation of Theorem \ref{thm2} can be obtained by considering (\ref{i2s}) under time-reversal. More specifically, we compare the state-equations for dynamical systems with transfer functions $V(\la)=D+C\la(I_n-\la A)^{-1}B$ and $V(\la)^*=D^*+B^*(\la I_n-A^*)^{-1}C^*$ given below:
\begin{eqnarray}\nonumber
\hspace*{-5pt}x_k&=&Ax_{k-1}+Bu_k\\
\hspace*{-5pt}y_k&=&Cx_{k-1}+Du_k,\; k=\ldots,-1,0,1,\ldots
\label{forward}
\end{eqnarray}
and
\begin{eqnarray}\nonumber
\hspace*{-19pt}x_{k-1}&=&A^*x_k+C^*y_k\\
\hspace*{-19pt}u_k&=&B^*x_k+D^*y_k,\; k=\ldots,1,0,-1,\ldots.
\label{backward}
\end{eqnarray}
Both are interpreted as stable linear dynamical systems but with opposite time-arrows.
Since $V(\la)V(\la)^*=I_m$, the input to one of the two corresponds to the output of the other, and
(\ref{spectralmeasures}) relates the spectral measure $d\mu$ of $\{u_k\}$ to the spectral measure $d\mu_r$ of $\{y_k\}$. The state-covariance for both system is the same when the first is driven by $\{u_k\}$ and the second by $\{y_k\}$, respectively. Thus, if $\bR=E\{x_kx_k^*\}$, Theorem \ref{thm1} applied to (\ref{forward})  leads to (\ref{rankcondition}a-b) while, applied to (\ref{backward}), leads to (\ref{rankcondition}c-d). The spectral measures of the respective inputs $\{u_k\}$ and $\{y_k\}$ relate as in (\ref{spectralmeasures}). 

\begin{proof}{\em [Theorem \ref{thm2}]} Follows readily from the above arguments. More precisely, $\bR$ is a state-covariance of (\ref{i2s}) for a suitable stationary input process $\{u_k\}$ if and only if it is also a state-covariance of
\[
x_{\ell+1}=A^*x_\ell+C^*y_\ell
\]
for a suitable stationary input process $\{y_\ell,\; \ell\in\mZ\}$. Then applying Theorem \ref{thm1} we draw the required conclusion.
\end{proof}

An analogous dual interpolation problem ensues. To avoid repeat of the development in \cite{acmatrix1,acmatrix2}, we may simply rewrite (\ref{Sigmar})
as
\begin{equation}\nonumber\label{Sigmar2}
\bR^\prime =\int_0^{2\pi}
\left(\left( G_r(e^{j\theta})^*\right)^\prime \frac{\left(d\mu_r(\theta)\right)^\prime}{2\pi} G_r(e^{j\theta})^\prime\right)
\end{equation}
where now the left integration kernel is
\[
\left( G_r(\la)^*\right)^\prime =(\la^{-1}I_n-A^\prime )^{-1}C^\prime.
\]
Note that $\bR^\prime=\bar{\bR}\neq \bR$ in general, since $\bR$ is Hermitian but may not be symmetric---where bar ($\bar{\;\;\;}$) denotes complex-conjugation.
Trading a factor $\la$ between the left integration kernel and its para-hermitian conjugate on the right we obtain that
\begin{equation}\nonumber\label{Sigmar2trade}
\hspace*{-7pt}\bR^\prime =\int_0^{2\pi}
(I_n-e^{j\theta} A^\prime )^{-1}C^\prime \frac{\left(d\mu_r(\theta)\right)^\prime}{2\pi}C(I_n-e^{-j\theta}A)^{-1}
\end{equation}
leading to the analytic interpolation problem of seeking an $\cF$-function of the form
\[
L^\prime(I_n-\la A^\prime)^{-1}C^\prime + Q(\la) V(\la)^\prime.
\]
Transposing once more we may define
\[
F_r(\la)=C(I_n-\la A)^{-1} L + V(\la)Q(\la)
\]
and draw the following conclusions.\\

\begin{thm}{\sf Let $V(\la)=D+C\la (I_n-\la A)^{-1}B$ be an $m\times m$ inner function with $(A,C)$ observable and $(A,B)$ reachable. If $L\in\mC^{n\times m}$ and $\bR$ the solution to (\ref{rankcondition}d), then there exists a solution $H$ to equation (\ref{rankcondition}b). Conversely, if $H\in\mC^{m\times n}$ and $\bR$ the solution to (\ref{rankcondition}b), then there exists a solution $L$ to equation (\ref{rankcondition}d). With $\bR,L,H$ related via (\ref{rankcondition}b) and (\ref{rankcondition}d), the following are equivalent:
\begin{eqnarray}
(\ref{thm3conditions}a)&&\bR\geq 0, \nonumber \\
(\ref{thm3conditions}b)&& \exists F(\la)\in\cF\;: \nonumber\\
&&F(\la)=H(I-\la A)^{-1}B+Q(\la)V(\la),\nonumber\\
&&\mbox{with }Q(\la)\mbox{ analytic in }\mD,\nonumber\\
(\ref{thm3conditions}c)&& \exists F_r(\la)\in\cF\;: \nonumber\\
&&F_r(\la)=C(I-\la A)^{-1}L+V(\la)Q_r(\la),
\nonumber\\
&&\mbox{with }Q_r(\la)\mbox{ analytic in }\mD.
\label{thm3conditions}
\end{eqnarray}
}\end{thm}
\vspace*{.1in}

\begin{proof} Begin with $L\in\mC^{n\times m}$ and $\bR$ the solution to (\ref{rankcondition}d). If $\bR\geq 0$ then $\bR$ is a state-covariance to (\ref{i2s}) according to Theorem \ref{thm2} and hence, there exists a solution $H$ to equation (\ref{rankcondition}b). To argue the case where $\bR$ may not be nonnegative definite necessarily, consider 
without loss of generality condition (\ref{standing1}e) valid and that
\[
U=\left[\begin{matrix}A &B\\C&D\end{matrix}\right]
\]
is unitary. Then, $I_n-A^*A=C^*C$ and $I_n-AA^*=BB^*$. If $\bR$ is the solution to (\ref{rankcondition}d) for a given $L$,
then $\bR_\epsilon:=\bR+\epsilon I_n$ is the solution of the same equation when $L$ is replaced by $L_\epsilon:=L+\frac{\epsilon}{2}C^*$. We can always choose $\epsilon$ so that $\bR_\epsilon>0$ and then deduce that there exists a solution $H_\epsilon$ to
\[
\bR_\epsilon-A\bR_\epsilon A^*=B H_\epsilon+H_\epsilon^* B^*.
\]
Since $I_n=AA^*+BB^*$, $H:=H_\epsilon-\frac{\epsilon}{2}B^*$ now satisfies (\ref{rankcondition}b). The converse proceeds in the same way.

The equivalence of (\ref{thm3conditions}a) and (\ref{thm3conditions}b) follows as in \cite{acmatrix1}.
If $\bR\geq 0$, then (\ref{thm3conditions}b) follows from \cite[Theorem 2]{acmatrix1}.
Conversely, if (\ref{thm3conditions}b) holds, then $\cW$ (defined in (\ref{cW}) satisfies (\ref{Wequation}) leading to $\bR$ being its Hermitian part. Since $F(\la)\in\cF$, the Hermitian part of multiplication by $F(\la)^*$ is nonnegative, and hence it remains so when restricted to the subspace $\cK$.

The dual statement (\ref{thm3conditions}c) follows in an analogous manner.
\end{proof}

\begin{remark}\label{smoothing} If $V(\la)=V_1(\la)V_2(\la)$ is a factorization of $V(\la)$ into a product of inner factors, then it can similarly be shown that the conditions (\ref{thm3conditions}) of the theorem are equivalent to the solvability of a bi-tangential Carath\'{e}odory-Fej\'{e}r interpolation problem of seeking an $F_o(\la)\in\cF$ where $F_o(\la)=H_o(I-\la A)^{-1}L_o+V_2(\la)Q(\la)V_1(\la)$ for suitable $H_o,L_o$. (The $H_o,L_o$ can be computed from e.g., $H,B$ by setting $F_o(\la)$ as the analytic part of $V_2(\la)F(\la)V_2(\la)^*$ and $F(\la)$ as in (\ref{thm3conditions}b).) $\Box$
\end{remark}

\section{\bf Optimal prediction \& postdiction errors}\label{predictionpostdiction}

A spectral distribution $\mu\in\mM$ induces
a Gram matricial structure on the space
of $p\times m$ matrix-valued functions on the circle
(see \cite[pages 353, 361]{masani}) via
\begin{eqnarray} \label{starintegral}
\hspace*{-21pt}\langle a(\la),b(\la)\rangle_{d\mu}&\hspace*{-9pt}:=&
\hspace*{-7pt}\int_0^{2\pi} b(e^{j\theta}) \frac{d\mu(\theta)}{2\pi}
a(e^{j\theta})^*\\
&\hspace*{-15pt}=& \hspace*{-15pt}\cE\{ (\sum_{\ell}b_\ell u_{k-\ell}) (\sum_{\ell}a_\ell u_{k-\ell})^*\}
\label{covinnerproduct}
\end{eqnarray}
where $a_\ell,b_\ell$ are the Laurent coeffients of $a(\la)$, $b(\la)$, respectively.
The correspondence
\begin{equation}\label{correspondence}
\sum_{\ell}a_\ell \la^\ell \mapsto \sum_{\ell}a_\ell u_{k-\ell},
\end{equation}
between functions on the unit circle (taking $\la=e^{j\theta}$) and
linear combinations of the random vectors $\{u_k\}$, leaves the respective Gram-matricial inner products in agreement and establishes a natural isomorphism between $\cL_2(\partial \mD;d\mu)$ and the space spanned by (the closure of) linear combination $\{u_k\}$
(see Masani \cite[Sections 5, 6]{masani}, cf.\ \cite{acmatrix2}).

Any
matrix-valued function
\begin{equation}\label{function_h}
h(\la) = \sum_{\ell=0}^\infty h_\ell \la^\ell
\end{equation}
with entries in ${\cal H}_2$ and
\begin{equation}\label{constraint}
h(0)=h_0=I_m
\end{equation}
corresponds via (\ref{correspondence}) to
\[
h(\la) \mapsto u_k-\hat{u}_{k|{\rm past}}
\]
which is interpreted as a ``one-step-ahead prediction error''.
Likewise, if the entries of
\begin{equation}\label{function_hr}
h(\la) = \sum_{\ell=0}^{-\infty} h_\ell \la^\ell
\end{equation}
live in $\la\cH_2^\perp$ and $h(0)=I_m$,
\[
h(\la) \mapsto u_k-\hat{u}_{k|{\rm future}}
\]
corresponds to ``one-step-ahead postdiction error'', i.e., using ``future'' observations only to determine the ``present''. Occasionally we may refer to these for emphasis as prediction forward, and backwards in time, respectively. Either way, the ``estimator'', which may not be optimal in any particular
way, is the respective linear combination of values of $u_k$ for $k\gtrless 0$:
\[
\hat{u}_{k|{\rm observation\;range}}:=-\sum_{\ell\gtrless 0} h_\ell u_{k-\ell}.
\]
(When the values extend in both directions it is a case of smoothing and is needed to interpret the $\cF$-function
in Remark \ref{smoothing}---this will be developed in a forthcoming report.)

We first discuss prediction in the forward direction. Throughout we consider as data the covariance matrix $\bR$ and the filter parameters. We assume that $d\mu\in\mM_\bR$ but otherwise unkown.
Because $d\mu$ is not known outside $\cK$, it can be shown that the min-max problem of identifying the forward prediction error with the least variance over all $d\mu\in\mM_\bR$ has a solution which lies in $\cK$. To this end we seek an element in $\cK^m$,
i.e., an $m\times n$ matrix-valued function
\[
\Cp G(\la) \mbox{ with } \Cp\in\mC^{m\times n}
\]
with rows in $\cK$, having least variance
\begin{eqnarray*}
\langle \Cp G(\la),\Cp G(\la)\rangle_{d\mu}
&=& \Cp \langle G(\la),G(\la)\rangle_{d\mu} \Cp^*\\
&=& \Cp \bR \Cp^*,
\end{eqnarray*}
and subject to the constraint (\ref{constraint}) which becomes
\begin{equation}\label{GammaBI}
\Cp B=I.
\end{equation}
Existence and characterization of minimizing matrices $\Cp$ is discussed next.

Nonnegative definiteness of the difference $\Omega_1-\Omega_2\geq 0$ between two elements
$\Omega_i\in\mH_m$ ($i=1,2$) defines a partial order $\Omega_1\geq \Omega_2$ in $\mH_m$. An $\mH_m$-valued function on a linear space
is said to be {\em $\mH_m$-convex} iff
\begin{eqnarray*}
&&f(\alpha \Cp_1+(1-\alpha)\Cp_2)\leq \alpha f(\Cp_1)+(1-\alpha)f(\Cp_2),\\
&&\mbox{for }\alpha\in[0,1].
\end{eqnarray*}
It is rather straightforward to check that if $\bR\geq 0$, then
the quadratic
\begin{equation}\label{f}
q_\bR\;:\;\mC^{m\times n}\to \mH_m \;:\; \Cp \mapsto \Omega=\Cp \bR \Cp^*,
\end{equation}
is in fact $\mH_m$-convex.
This basic fact ensures existence of $\mH_m$-minimizers satisfying (\ref{GammaBI}) in the proposition given below.
Note that the statements (ii) and (iii) of the proposition are rephrased in alternative ways (e.g., (ii-a), etc.) in order to highlight an apparent symmetry when expressed in terms of directed gaps $\vec{\delta}$ (defined in the statement of the proposition) between
the null space
\[
\cN(\bR):=\{ x\in\mC^{n\times 1}\;:\;\bR x=0_n\}
\]
of $\bR$ and the range
\[
\cR(B):=\{ x\in\mC^{n\times 1}\;:\; x=Bv \mbox{ for }v\in\mC^{m\times 1}\}
\]
of $B$---the gap metric represents an angular distance between subspaces and is a standard tool in perturbation theory of linear operators (see \cite{Kato}) and in robust control (e.g., see \cite{mopus}).\\

\begin{prop}\label{lemma0}
{\sf 
Let $B\in\mC^{n\times m}$ having rank $m$, and let $\bR\in\mH_n$ with $\bR\geq 0$.
The following hold:
\begin{itemize}
\item[(i)] There exists an $\mH_m$-minimizer of $q_\bR$ satisfying (\ref{GammaBI}).
\item[(ii)] The minimizer is unique if and only if
\[
\rank(\left[\begin{array}{cc}\bR& B\end{array}\right])=n.
\]
\item[(ii-a)] The minimizer is unique if and only if
\[
\vec{\delta}(\cN(\bR),\cR(B)):=\| \bPi_{\cR(B)^\perp}|_{\cN(\bR)}\|<1.
\]
\item[(iii)] The $\mH_m$-minimal value for $q_\bR$ is $O_m$ if and only if
\[
B^*\bPi_{\cN(\bR)}B \mbox{ is invertible}.
\]
\item[(iii-a)] The $\mH_m$-minimal value for $q_\bR$ is $O_m$ if and only if
\[
\vec{\delta}(\cR(B),\cN(\bR))):=\| \bPi_{\cN(\bR)^\perp}|_{\cR(B)}\|<1.
\]
\item[(iv)] If $\rank(\bR)=n$, then the $\mH_m$-minimal value of $q_\bR$ is
\[
\Omega:=(B^*\bR^{-1} B)^{-1}>0
\]
and a minimizer (unique by (ii)) is
\[
\Cp=(B^*\bR^{-1} B)^{-1}B^*\bR^{-1}.
\]
\item[(v)] If the $\mH_m$-minimal value of $q_\bR$ is
\[
\Omega=O_m,
\]
then a minimizer is given by
\begin{equation}\label{Gamma0}
\Cp=(B^*\bPi_{\cN(\bR)} B)^{-1}B^*\bPi_{\cN(\bR)}.
\end{equation}
\item[(vi)] In general, when $\bR$ is singular, the $\mH_m$-minimal value for $q_\bR$ is
\begin{equation}\label{Omegageneral}
\Omega=(B_1^*\bR^\sharp B_1)^\sharp
\end{equation}
and a minimizer is given by
\begin{eqnarray}\label{Gammageneral}
\Cp&=&(B^*\bPi_{\cN(\bR)} B)^\sharp B^*\bPi_{\cN(\bR)}\\\nonumber
&&\hspace*{-2.5cm}+(B_1^*\bR^\sharp B_1)^\sharp B_1^*\bR^\sharp(I_n-B(B^*\bPi_{\cN(\bR)} B)^\sharp B^*\bPi_{\cN(\bR)})
\end{eqnarray}
where $\bR^\sharp$ denotes the Moore-Penrose pseudo-inverse of $\bR$,
\begin{eqnarray*}
B_1&:=&B\bPi_{\cN_o}, \mbox{ and}\\ \cN_o&:=&\cN(B^*\bPi_{\cN(\bR)}B).
\end{eqnarray*}
Alternatively,
\begin{eqnarray}\nonumber
\Omega&=&\lim_{\epsilon\to 0}\Omega_\epsilon\\
\Gamma&=&\lim_{\epsilon\to 0}\Gamma_\epsilon \label{alt1}
\end{eqnarray}
where
\begin{eqnarray}\nonumber
\Omega_\epsilon&:=&(B^*\bR_\epsilon^{-1}B)^{-1}\\\nonumber
\Gamma_\epsilon &:=&\Omega_\epsilon^{-1}B^* \bR_\epsilon^{-1}\\
\bR_\epsilon&:=&\bR+\epsilon \bPi_{\cN(\bR)} \label{alt2}\end{eqnarray}
and $\epsilon>0$.
\end{itemize}
}
\end{prop}
\vspace*{.1in}

\begin{proof}
{\bf Claim (i):} Since $q_\bR(\Cp)\geq 0$, then for any $\alpha\in[0,1]$ and any $\Cp_1,\Cp_2$
\begin{eqnarray*}
&&\alpha(1-\alpha)q_\bR(\Cp_1-\Cp_2)\geq 0\\[.05in]
&\Leftrightarrow& \alpha(1-\alpha)\left(
\Cp_1\bR\Cp_1^*+\Cp_2\bR\Cp_2^*\right. \\ &&\left.-\Cp_1\bR\Cp_2^*-\Cp_2\bR\Cp_1^*
\right)\geq 0\\[.05in]
&\Leftrightarrow&
\alpha\Cp_1\bR\Cp_1^*
+(1-\alpha)\Cp_2\bR\Cp_2^*\\
&&
-\alpha^2\Cp_1\bR\Cp_1^*
-(1-\alpha)^2\Cp_2\bR\Cp_2^*
\\ &&-\alpha(1-\alpha)\Cp_1\bR\Cp_2^*-\alpha(1-\alpha)\Cp_2\bR\Cp_1^* \geq 0\\[.05in]
&\Leftrightarrow&
\alpha q_\bR(\Cp_1)+(1-\alpha) q_\bR(\Cp_2)\\ &&\geq  q_\bR(\alpha\Cp_1+(1-\alpha)\Cp_2).
\end{eqnarray*}
This proves $\mH_m$-convexity of $q_\bR$. It is also clear that $q_\bR$ is bounded below by $O_m$. However, $q_\bR$ is not necessarily radially unbounded when $\bR$ is singular. Hence, we need to consider components of $\Gamma$ which lie in $\cN(\bR)$. Any $\Gamma$ satisfying (\ref{GammaBI}) is of the form
\[
\Gamma=\Gamma_0+X M
\]
where $\Gamma_0$ is a particular solution of (\ref{GammaBI}) (e.g., $\Gamma_0=(B^*B)^{-1}B^*$), the rows of $M\in\mC^{(n-m)\times n}$ span the left null space of $B$, and $X$ is an arbitrary element of $\mC^{m\times(n-m)}$. Substituting into $q_\bR$ we obtain
\newcommand{\bQ}{{\bf Q}}
\newcommand{\bL}{{\bf L}}
\newcommand{\bC}{{\bf C}}
\begin{equation}\label{generalquadratic}
q_\bR(\Gamma_0+XM)= X\bQ X^* + X\bL + \bL^* X^* + \bC,
\end{equation}
with $\bQ=M\bR M^*$, $\bL=M\bR \Gamma_0^*$, and $\bC=\Gamma_0 \bR \Gamma_0^*$,
which is also $\mH_m$-convex in $X$. Since it is bounded below by $O_m$, it follows that the null space of $\bQ$ is contained in the null space of $\bL^*$. Expressing the entries in (\ref{generalquadratic}) with respect to the decomposition $\mC^n=\cR(\bQ)\oplus \cN(\bQ)$, we may write $q_\bR(\Gamma_0+XM)$ in the form
\begin{eqnarray}
&&\left[ \begin{matrix}X_1&X_2\end{matrix}\right] \nonumber
\left[\begin{array}{cc}\bQ_1& O_{(n-m)\times m}\\
O_{m\times(n-m)}&O_{m}\end{array}\right]
\left[ \begin{matrix}X_1^*\\X_2^*\end{matrix}\right]\\
&&\hspace*{-17pt}+ \left[ \begin{matrix}X_1&X_2\end{matrix}\right]\left[ \begin{matrix}\bL_1\\O_m \end{matrix}\right] + \left[ \begin{matrix}\bL_1^*&O_m\end{matrix}\right] \left[ \begin{matrix}X_1^*\\X_2^*\end{matrix}\right] + \bC,\label{generalquadratic2}
\end{eqnarray}
where $\bQ_1>0$. This expression is radially unbounded in $X_1$
and hence, a minimizer exists (taking any bounded value for $X_2$, e.g., $X_2=O_m$). This proves (i).

{\bf Claim (ii):} Because all minimizers satisfy (\ref{GammaBI}), any two of them differ by some matrix, say $\Delta$ such that $\Delta B=O_m$. Hence if $q_\bR(\Cp)=q_\bR(\Cp+\Delta)$, then $q_\bR(\Cp)=q_\bR(\Cp+\epsilon\Delta)$, for $\epsilon\in[0,1]$. This is due to the $\mH_m$-convexity of $q_\bR$. Therefore
\[
\epsilon^2 \Delta\bR\Delta^*+ \epsilon \Cp \bR \Delta^* + \epsilon \Delta \bR \Cp^*=O_m,
\]
identically for all $\epsilon\in[0,1]$. It follows that there is more than one minimizer if and only if
there exists a common left null vector  for both $B$ and $\bR$ (which serves as a nonzero row of $\Delta$, so that $\Delta\neq O_{m\times n}$). This proves (ii).

{\bf Claim (ii-a):} For the definition of the directed gap in (ii-a) cf.\  \cite{Kato}.
The claim that (ii-a) is equivalent to (ii) is standard.
Since $\cR(B)^\perp$ coincides with $\cN(B^*)$, a common element between $\cN(B^*)$ and $\cN(\bR)$ would lead to $\| \bPi_{\cR(B)^\perp}|_{\cN(\bR)}\|=1$. Since we are dealing with finite-dimensiional spaces the converse is immediate---a common vector is the only way the norm can be equal to one in this case. The rank condition in (ii) is obviously equivalent to $\cN(B^*)\cap\cN(\bR)=\{0\}$.

{\bf Claim (iii) and claim (v):} We now argue claim (iii) together with
claim (v). If $B^*\bPi_{\cN(\bR)}B$ is invertible and $\Cp$ as in (\ref{Gamma0}), then $\Cp B=I_m$ and $q_\bR(\Cp)=O_m$.
To show the converse assume that $\exists \Cp$ such that $q_\bR(\Cp)=O_m$ as well as $\Cp B=I_m$. Then the columns of $\Cp^*$ belong to $\cN(\bR)$ and $\Cp \bPi_{\cN(\bR)}=\Cp$. Therefore
$\Cp \bPi_{\cN(\bR)}B=I_m$ and $\bPi_{\cN(\bR)}B$ has rank $m$. Consequently, $B^*\bPi_{\cN(\bR)}B$ is invertible.

{\bf Claim (iii-a):} The equivalence of (iii-a) and (iii) is standard. Condition (iii-a) is equivalent to stating that
$\bPi_{\cN(\bR)}|_{\cR(B)}$ has rank $m$. But $\bPi_{\cN(\bR)}|_{\cR(B)}+\bPi_{\cR(\bR)}|_{\cR(B)}=I_{\cR(B)}$ (where $I_{\cR(B)}$ denotes the identity operator on $\cR(B)$).
Because $\bR\in\mH_n$, $\cN(\bR)^\perp=\cR(\bR)$, and condition (iii) follows.

{\bf Claim (iv):} Assume that $\bR$ is positive definite and $\Gamma,\,\Omega$ as in (iv). Then
$\Gamma B=I_m$ and $q_\bR(\Gamma)=\Omega>0$.
For any $X\in\mC^{m\times n}$ such that $XB=O_m$, it can be readily seen that $q_\bR(\Gamma+X)=\Omega+X\bR X^*>\Omega$.
Hence, the minimizer and minimal value are as claimed. This proves (iv).

{\bf Claim (vi):} Denote
\[
\cR_o:=\cR(B^*\bPi_{\cN(\bR)}B)
\]
and recall that, for any matrix $M$, the orthogonal projection onto its range can be obtained via
\[
\bPi_{\cR(M)}=M^\sharp M.
\]
We verify by direct substitution that
\begin{eqnarray}
\hspace*{-2cm}\Gamma B&=&(B^*\bPi_{\cN(\bR)} B)^\sharp B^*\bPi_{\cN(\bR)}B\nonumber
+(B_1^*\bR^\sharp B_1)^\sharp(B_1^*\bR^\sharp B_1)\\
&=&\bPi_{\cR_o}+\bPi_{\cN_o} \label{step}\\&=& I_n,\nonumber
\end{eqnarray}
and that
\begin{eqnarray*}
\Gamma R\Gamma^*&=& (B_1^*\bR^\sharp B_1)^\sharp(B_1^*\bR^\sharp B_1)(B_1^*\bR^\sharp B_1)^\sharp\\
&=&\Omega
\end{eqnarray*}
as given in (vi). Step (\ref{step}) needs the fact that
\[\cR(B_1^*\bR^\sharp B_1)=\cN_o.
\]
We can show this as follows. Clearly $\cR(B_1^*\bR^\sharp B_1)\subseteq \cN_o$ since $B_1^*=\bPi_{\cN_o}B^*$.
To establish equality we need to show that there exists no $x\in\cN_o$ other than $0$  such that $R^\sharp Bx=0$, i.e., that
\begin{equation}\label{intersection}
\cN(B^*R^\sharp B)\cap \cN_o=\{0\}.
\end{equation}
But
\[
\cN(B^*R^\sharp B)=\cN(B^*\bPi_{\cR(\bR)} B)
\]
and
\begin{eqnarray*}
B^*\bPi_{\cR(\bR)} B+B^*\bPi_{\cN(\bR)} B&=&B^*(\bPi_{\cR(\bR)}+\bPi_{\cN(\bR)})B\\
&=&B^*B
\end{eqnarray*}
is invertible. Hence (\ref{intersection}) holds and so does (\ref{step}).

We finally need to show that the value for $\Omega$ is an $\mH_m$-minimum of $q_\bR$ subject to (\ref{GammaBI}). For any $X$ such that $XB=O_m$ it also holds that $XB_1=O_m$. We can verify by direct substitution that
\[
q_\bR(\Gamma+X)=\Omega+X\bR X^*,
\]
which proves that $\Omega,\,\Gamma$ as given represent the minimum and minimizer, respectively.

We argue the validity of the alternative set of expressions (\ref{alt1}-\ref{alt2}) as follows.
For any $\epsilon\in(0,1]$,
\begin{equation}\label{Sigmaepsilon}
\bR_\epsilon=\bR+\epsilon \bPi_{\cN(\bR)}
\end{equation}
is positive definite with
\[\bR_\epsilon^{-1}=\bR^\sharp +\epsilon^{-1}\bPi_{\cN(\bR)}
\]
as its inverse. We can now apply (iv) to argue that $\Omega_\epsilon,\,\Gamma_\epsilon$ are the minimal value and minimizer of $q_{\bR_\epsilon}$ subject to $\Gamma_\epsilon B=I_m$, as before. It follows that their limits satisfy $\Gamma B=I_m$ and
$\Gamma \bR\Gamma^*=\Omega$. 
Then $\Omega$ is indeed the $\mH_m$-minimal value of $q_\bR$ (cf.\ (i)) by continuity of $q_{\bR}$ on $\bR$.

It is straightforward (but a bit cumbersome to typeset) to use the limits (\ref{alt1}) and verify (\ref{Omegageneral}-\ref{Gammageneral}).
To pursue this, express $B^*\bR_\epsilon^{-1}B$ as a $2\times 2$ matrix with respect to the decomposition
\[ \mC^m=\cN_o\oplus \cR_o.
\]
The (1,1) entry, $\alpha:=\bPi_{\cN_o} B^*\bR^\sharp B|_{\cN_o}$ is invertible and so is the (2,2) entry
\[
\bPi_{\cR_o} B^*\bR^\sharp B|_{\cR_o}+ \epsilon^{-1}B^*\bPi_{\cN(\bR)}B =:\gamma+\epsilon^{-1}\delta
\]
where $\gamma,\delta$ are defined to represent the respective terms.
The (2,2) entry is the only one involving the parameter $\epsilon$. Then, the inverse of $B^*\bR_\epsilon^{-1}B$
becomes
\[
\Omega_\epsilon =\left[\begin{array}{cc} \alpha^{-1}+o(\epsilon) &-\alpha^{-1}\beta\delta^{-1}\epsilon +o(\epsilon^2)\\
-\delta^{-1}\beta^*\alpha\epsilon +o(\epsilon^2) & \epsilon\delta^{-1}+o(\epsilon^2)\end{array}\right]
\]
with $\beta:=\bPi_{\cN_o} B^*\bR^\sharp B|_{\cR_o}$. The limit gives the correct expression for $\Omega$. The limit of $\Gamma_\epsilon$ as $\epsilon\to 0$ can be carried out similarly.
\end{proof}

\begin{remark}It should be noted that $\bR$
is not required to have the structure of a state-covariance of a reachable pair $(A,B)$ (cf.\ Theorem \ref{thm1}) since
the matrix $A$ does not enter at all in the statement of Proposition \ref{lemma0} . However, if this is the case (see Proposition \ref{Omegasingular} below) and $\bR$ is a singular \underline{state-covariance}, then $\Omega$ is singular as well---a converse to the first part of statement (iv).
$\Box$ \end{remark}

For prediction backwards in time, the postdiction error
\begin{equation}\label{postdiction}
u_k-\hat{u}_{k|future} = u_k-\sum_{\ell=1}^\infty h_\ell u_{k+\ell}
\end{equation}
corresponds to an element
\begin{eqnarray*}
&&zL^*G_r(z)=L^*(I_n-z^{-1}A^*)^{-1}C^*\in z\cK_r\\
&& \mbox{with }L\in\mC^{n\times m}.
\end{eqnarray*}
The constraint arising from the the identity in front of $u_k$ in (\ref{postdiction}), translates into
\[
L^*C^*=I_m
\]
while the variance of the postdiction error becomes
\[
L^*\bR L.
\]
Proposition \ref{lemma0} applies verbatim and yields that:
\begin{itemize}
\item[(i')] there exists an $\mH_m$-minimal postdiction error.
\item[(ii')] The minimizer is unique if and only if
\[\rank(\left[\begin{array}{c}\bR\\ C\end{array}\right])=n. \]
\item[(iii')] The variance of optimal postdiction error is equal to $O_m$ if and only if
\[C\bPi_{\cN(\bR)}C^* \mbox{ is invertible}.\]
\item[(iv')] If $\rank(\bR)=n$, then the variance of the optimal postdiction error is (strictly) positive definite and
the unique minimizer is
\[\Gamma_r=\bR^{-1}C^*(C\bR^{-1} C^*)^{-1}.
\]
\item[(v')] If the variance of the optimal postdiction error is equal to $O_m$, then a (non-unique) minimizer is
\begin{equation}\label{Gamma0prime}
\Gamma_r=\bPi_{\cN(\bR)}C^*(C\bPi_{\cN(\bR)} C^*)^{-1}.
\end{equation}
\end{itemize}
Similarly, the analog of (vi) holds as well.


\begin{remark} It is interesting to point out that the square-roots of the variances of prediction and postdiction errors $(B^*\bR^{-1}B)^{-1/2}$ and $(C\bR^{-1}C^*)^{-1/2}$ appear as left and right radii, respectively, in a Schur parametrization of the elements of $\mM_\bR$ in \cite{acmatrix1} (cf.\ \cite[Remark 2]{acmatrix2}) and that, in view of the above, if one is zero so is the other. $\Box$ \end{remark}

\section{\bf When $\mM_\bR$ contains a single element}\label{sec:singleton}

We now focus on the case where $\mM_\bR$ consists of a single element, we
analyze the nature of this unique power spectrum, and study ways to decompose $\bR$ into a sum of two non-negative definite matrices, one of which has this property and another which may be interpreted as corresponding to noise.  Conditions for $\mM_\bR$ to be a singleton are stated next.\\

\begin{thm}\label{mainthm1}{\sf
Let $A,B$ satisfy (\ref{standing1}) and $\bR\geq 0$ for which (\ref{rankcondition}) holds.
Then, the set $\mM_\bR$ is a singleton if and only if the following equivalent conditions
hold:
\begin{eqnarray}
(\ref{singleton}a) && \vec{\delta}(\cR(B),\cN(\bR))<1,\nonumber\\
(\ref{singleton}b) && B^*\bPi_{\cN(\bR)}B \mbox{ is invertible.}\label{singleton}
\end{eqnarray}
If $(C,D)$ are selected so that $V(\la)$ in (\ref{V})  is inner, the above conditions are also equivalent to:
\begin{eqnarray}
(\ref{singleton}c) && \vec{\delta}(\cR(C^*),\cN(\bR))<1,\nonumber\\
(\ref{singleton}d) && C\bPi_{\cN(\bR)}C^* \mbox{ is invertible.}\nonumber
\end{eqnarray}
}\end{thm}
\vspace*{.1in}

\begin{proof}
As explained earlier, an element $d\mu\in\mM$
defines via (\ref{correspondingpr}) an $\cF$-function
$F(\la)=\cH[d\mu]$ which, in turn, defines a (possibly unbounded) non-negative operator on $\cH_2^{1 \times m}$ via
\[
x(\la) \mapsto \bPi_{\cH_2} x(\la)F(\la)^*.
\]
Conversely, this operator defines uniquely the function $F(\la)\in\cF$ as well as the corresponding measure $d\mu\in\mM$ (except of course for a skew-Hermitian constant in $F(\la)$ and an additive constant in $\mu$).
The restriction onto $\cK$,
\begin{eqnarray*}
\cW : \cK\to \cK&:& x(\la) \mapsto \bPi_{\cK} x(\la)F(\la)^*
\end{eqnarray*}
corresponds to ``one half'' of $\bR$ as in (\ref{Wequation}), and is specified by $\bR$ (modulo a skew-Hermitian part). We proceed to show recursively that there exists a unique extension of $\cW$ to a non-negative operator on
\[
\cK_\ell := \cH_2^{1 \times m}\ominus  \cH_2^{1 \times m}\la^\ell V(\la)
\]
for $\ell=1,\,2,\ldots$, and hence, to a non-negative operator on $\cH_2^{1 \times m}$.

Consider the representation
\[
\la V(\la)=D_1+C_1\la (I-\la A_1)^{-1}B_1
\]
with
\begin{eqnarray} \label{notationABCD}
{\cal A}_1&=& \left[\begin{array}{cc} A & 0\\
                               C & 0\end{array}\right],\\\nonumber
{\cal B}_1&=& \left[\begin{array}{c}  B \\
                               D \end{array}\right], \\\nonumber
{\cal C}_1&=& \left[\begin{array}{cc} 0 & I\end{array}\right],\\\nonumber
{\cal D}_1&=& 0,
\end{eqnarray}
and
\[
\bR_1:=\left[\begin{array}{cc}\bR &\bR_{12}\\\bR_{12}^* & \bR_{22}\end{array}\right]
\]
the (non-negative) Hermitian part of an extension of $\cW$ into $\cK_1$. Then, from Theorem \ref{thm1},
\begin{equation}\label{extension}
\bR_1-{\cal A}_1\bR_1 {\cal A}_1^*={\cal B}_1 {\cal H}_1 + {\cal H}_1^*{\cal B}^*
\end{equation}
where
\begin{equation}\label{cH1}
{\cal H}_1=\left[\begin{array}{cc}H & H_1 \end{array}\right].
\end{equation}

Let us first assume that (\ref{singleton}a) holds (and hence, from Proposition \ref{lemma0}, that (\ref{singleton}b-d) hold as well).
Then $\Gamma \bR=O_{m\times n}$ with $\Gamma$ as in (\ref{Gamma0}) satisfying $\Gamma B=I_m$. Because, $\bR_1\geq 0$, it follows that
\[
\Gamma \bR_{12}=O_{m\times m},
\]
otherwise it would  be possible to render the quadratic form $\alpha\Gamma \bR_{12}+\bar{\alpha}\bR_{12}^*\Gamma^*+|\alpha|^2\bR_{22}$ indefinite with a suitable choice of $\alpha\in\mC$ which would contradict $\bR_1\geq 0$.
From (\ref{extension}) on the other hand, we have that
\[
\bR_{12}-A\bR C^*=B H_1+ H^* D^*.
\]
Multiplying on the left and the right by $\Gamma$ and $B$, respectively, we conclude that
\[ H_1= -\Gamma A\bR C^*-\Gamma H^* D^*
\]
is uniquely defined from the original data---hence, so is the ``one-step'' extension $\bR_1$ of $\bR$. It remains to show that the condition (\ref{singleton}a) is still valid for the new data, i.e., that
\[
{\cal B}^*_1\bPi_{\cN(\bR_1)}{\cal B}_1
\]
is also invertible. Since $\bR$ is Hermitian,
\[
\mC^n=\cN(\bR)\oplus \cR(\bR)
\]
is an orthogonal decomposition. Then, the null space of $\bR_1$ is the orthogonal direct sum of
\[
\{ \left(\begin{array}{c} \bPi_{\cN(\bR)}x\\ 0\end{array}\right)\;:\;x\in\mC^n \}
\]
and
\[
\{ \xi=\left(\begin{array}{c} \bPi_{\cR(\bR)}x\\ y\end{array}\right)\;:\; x\in\mC^n,\;y\in\mC^m,\; \bR_1\xi=0\}.
\]
Denote these two subspaces by $\cN_1$ and $\cN_2$, respectively. Then,
\[
\bPi_{\cN(\bR_1)}=\bPi_{\cN_1}+\bPi_{\cN_2}
\]
where
\[
\bPi_{\cN_1}=\bPi_{\cN(\bR)}\oplus O_{m\times m}.
\]
So, finally,
\[
{\cal B}^*_1\bPi_{\cN(\bR_1)}{\cal B}_1=B^*\bPi_{\cN(\bR)}B+{\cal B}^*_1\bPi_{\cN_2}{\cal B}_1 > 0,
\]
because $B^*\bPi_{\cN(\bR)}B$ is already positive definite. This completes the proof.
\end{proof}

The unique element in $\mM_\bR$ under the conditions of the theorem can be obtained, in principle, after extending $\bR$ recursively for $\ell=1,2,\ldots$ using (\ref{notationABCD}-\ref{cH1}). This specifies a non-negative operator on a dense subset of $\cH_2^{1\times m}$ which, in turn, specifies a corresponding positive real function $F(\la)$ and the measure can be obtained from the boundary limits of the real part of $F(\la)$ as a weak limit. However, an explicit expression for $F(\la)$ will also be given later on. Before we do this we explain some of the properties of this unique measure.

The following result states that $d\mu$ is a singular measure with at most $n=m$ points of increase, i.e., at most $n-m$ spectral lines whose directionality is encapsulated in suitably chosen unitary factors. The spectral lines are in fact at the zeros of certain matrix-valued functions, namely
\begin{equation}\label{Phi}
\Phi(\la):=\Gamma G(\la) =\Gamma(I_n-\la A)^{-1}B
\end{equation}
and $\Gamma$ as in Proposition \ref{lemma0},
which correspond to the optimal prediction error and represent the analog of the Szeg\"{o}-Geronimus orthogonal polynomials of the first kind, cf.\ \cite{acmatrix2}.

\begin{thm}\label{mainthm2}{\sf
Under the assumptions and conditions of Theorem \ref{mainthm1}, the unique element in $\mM_\bR$ is of the form
\[
d\mu(\theta)=\sum_{\ell=1}^q V_\ell \rho_\ell V_\ell^* d\mU(\theta-\theta_\ell)
\]
where
\[
\sum_1^q\rank(V_\ell)\leq n-m,
\]
$\theta_\ell\in[0,2\pi)$ for $\ell=1,\ldots,q$ differ from one another, $\mU(\theta-\theta_\ell)$ denotes a unit step at $\theta_\ell$, and $\rho_\ell>0$. The values $e^{j\theta_\ell}$ for $\ell=1,\ldots,q$ are the non-zero eigenvalues of the matrix
$(I_n-B\Cp)A$ with $\Cp$ as in (\ref{Gamma0}). The matrices $V_\ell$ are chosen so that
\[
\cR(V_\ell)=
\cN(B^*\bPi_{\cN(\bR)}(I_n-e^{j\theta_\ell}A)^{-1}B),
\]
and can be normalized to satisfy $V_\ell V_\ell^*=I$ as well as to make $\rho_\ell$ diagonal.
}\end{thm}
\vspace*{.1in}

\begin{proof}
Under the stated conditions, $\mM_\bR$ is a singleton from the previous theorem and
its unique element $d\mu$ satisfies
\begin{equation}\label{error0}
\int_0^{2\pi} \left(\Phi(e^{j\theta})d\mu(\theta)\Phi(e^{j\theta})^*\right)=
\Gamma \bR \Gamma^*=
O_{m\times m}
\end{equation}
with $\Gamma$ as in (\ref{Gamma0}).
It readily follows that $d\mu$ can have points of increase only at the finitely many points $\theta_\ell$, $\ell=1,\ldots,q$, where $\Phi(e^{j\theta})$ is singular. The ``zeros'' of $\Phi(\la)$ coincide with the ``poles'' of its inverse
\begin{equation}\label{Phiinv}
\Phi(\la)^{-1} = I_n-\Gamma A(I_n-\la A_o)^{-1}B
\end{equation}
where
\begin{equation}\label{Ao}
A_o=(I_n-B\Gamma)A.
\end{equation}
Since $A_o$ has already $m$ eigenvalues at the origin, the number of eigenvalues that it may have on the circle is at most $n-m$.
Thus
\[
d\mu(\theta)=\sum_{\ell=1}^q M_\ell d\mU(\theta-\theta_\ell)
\]
where $M_\ell\in\mH_m$, $M_\ell\geq 0$, and
\[
\Phi(e^{j\theta_\ell})M_\ell = O_{m\times m}.
\]
Exressing $M_\ell=V_\ell \rho_\ell V_\ell^*$ with $\rho_\ell,V_\ell$ as claimed is standard.
This completes the proof.
\end{proof}

Thus, $\mM_\bR$ being a singleton implies just as in the classical scalar case (e.g., \cite{StoicaMoses,GrenanderSzego}) that the underlying stochastic process is deterministic with finitely many complex exponential components. Subspace identification techniques  represent different ways to identify ``dominant ones'' and obtain the ``residue'' $\rho_\ell$ that corresponds to each of those modes  (see \cite{GrenanderSzego}, \cite{StoicaMoses}, \cite{sesha1,sesha2}). In the present multivariable setting, in order to do something analogous, we need an explicit expression for the corresponding positive real function. This is done in the next section.

\begin{remark}
A dual version of the representation in Theorem \ref{mainthm2} gives that
$\theta_\ell$ correspond to ``zeros'' on the circle of the optimal postdictor error
\[
L^*(\la I_n-A^*)^{-1}C^*.
\]
Similarly, the range $\cR(V_\ell)$, for $\ell=1,2,\ldots,q$, is contained in the correspond null space of the above postdiction error when evaluated at the corresponding zeros.
$\Box$ \end{remark}
\begin{remark}
The ``star'' of the optimal postiction error can also be interpreted as a ``right matricial orthogonal polynomial of the first kind''
\begin{equation}\label{Phir}
\Phi_r(\la)=C(I_n-\la A)^{-1}L.
\end{equation}
These matricial functions, i.e., $\Phi(\la)$ and $\Phi_r(\la)$, together with their counterparts of the ``second kind'' $\Psi(\la)$ and $\Psi_r(\la)$ that will be introduced in the next section, satisfy a number of interesting properties similar to those of the classical orthogonal polynomials \cite{Geronimus} (cf.\ \cite{dgk1,dgk2}). We plan to develop this subject in a separate future publication. $\Box$
\end{remark}

\section{\bf The ``central'' positive real function}\label{positiverealfunction}

With $\bR,A,B,H$ satisfying (\ref{rankcondition}b) in Theorem \ref{thm2}, we define
\newcommand{\Dpsi}{D_\Psi}
\begin{equation}
F_{ME}(\la):= \Phi(\la)^{-1}\Psi(\la) \label{FME}
\end{equation}
where $\Phi(\la)=\Gamma(I_n-\la A)^{-1}B$ as before, 
\begin{equation}\label{Psi}
\Psi(\la):=-\Gamma\la(I_n-\la A)^{-1}AH^* + \Dpsi,
\end{equation}
and
\begin{equation}\label{Dpsi}
\Dpsi:= -\Gamma(H^*B^*-\bR)B(B^*B)^{-1}.
\end{equation}
By eliminating the unobservable dynamics in the expression for $F_{ME}(\la)$ we obtain
\begin{eqnarray}
F_{ME}(\la)&=& \Dpsi+\la C_o(I_n-A_o)^{-1}B_o \label{reducedFME}
\end{eqnarray}
where
\begin{eqnarray}\nonumber
C_o&:=&-\Gamma A\\\nonumber
A_o&:=&(I_n-B\Gamma)A\\
B_o&:=&BD_\Psi+H^*.\label{reducedMEdata}
\end{eqnarray}

In case $\bR>0$, $F_{ME}(\la)$ is the positive-real functions which corresponds to the ``maximum entropy'' spectral measure $d\mu_{ME}(\theta)\in\mM_\bR$, i.e., the unique element of $\mM_\bR$ which maximizes the entropy functional
\[
{\mathbb{I}}(\mu):=\int_0^{2\pi} \log\det \left(\dot{\mu}(\theta)\right) d\theta.
\]
This element was identified in \cite{acmatrix2} as
\begin{equation}\label{mudot}
d\mu_{ME}=\Phi(e^{j\theta})^{-1}\Omega \left(\Phi(e^{j\theta})^{-1}\right)^*d\theta,
\end{equation}
without drawing the connection to (\ref{reducedFME}).
However, $F_{ME}(\la)$ in (\ref{FME}) is defined even when $\bR$ is singular, in which case the corresponding measure may have a singular part obtained as the weak radial limit of the Hermitian part of $F_{ME}(\la)$ 
\begin{equation}\label{radiallimits2}\nonumber
d\mu_{ME}(\theta) = \lim_{r\nearrow 1}{\Herm}\{F_{ME}(re^{j\theta})\}d\theta.
\end{equation}
The singular part, which corresponds to purely deterministic components in the underlying time series, relates to the residues of $F_{ME}(\la)$ at corresponding poles on the unit circle. This allows identifying spectral lines directly from $F_{ME}(\la)$. It should be emphasized that (\ref{mudot}) is no longer valid in the case of a singular $\bR$. We first establish the claim that $F_{ME}$ is positive real and that it is consistent with $\bR$.\\

\begin{thm}\label{thm6}{\sf
Let $\bR,A,B,H$ satisfy (\ref{rankcondition}b) of Theorem \ref{thm2}, $\Gamma$ given as in (\ref{Gammageneral}), and $F_{ME}(\la)$ given as in (\ref{FME}). Then
\begin{itemize}
\item[(i)] $F_{ME}(\la)$ satisfies (\ref{PRNehari}), and
\item[(ii)]  $F_{ME}(\la)\in\cF$.
\end{itemize}
}\end{thm}
\vspace*{.1in}

\begin{proof}
Condition (\ref{PRNehari}) is equivalent to
\[\Psi(\la)V(\la)^*-\Phi(\la)HG(\la)V(\la)^*=\Phi(\la)Q(\la).
\]
To show that this relationship holds for some $Q(\la)$ analytic in $\mD$, it suffices to show that all negative Fourier coefficients of
\begin{equation}\label{negativecoef}
\Psi(\la)V(\la)^*-\Phi(\la)HG(\la)V(\la)^*
\end{equation}
 vanish. By collecting positive and negative powers of $\la$ we can express
\begin{eqnarray*}
\Psi(\la)V(\la)^*&=&(D_\Psi B^*-\Gamma A W A^*)(\la I_n-A^*)^{-1}C^*\\ &&
 \hspace*{-2cm}+ D_\Psi D^*-\Gamma A (I_n-\la A)^{-1}(\la H^*D^*+WC^*),
\end{eqnarray*}
and similarly that
\begin{eqnarray*}\Phi(\la)HG(\la)V(\la)^*&=&\Gamma W^*(\la I_n-A^*)^{-1}C^*\\
&&+\Gamma (I_n-\la A)^{-1}WC^*
\end{eqnarray*}
where $W$ is given by (\ref{Wequation}).
Thus, negative powers of $\la$ in (\ref{negativecoef}) sum up into
\[\left(D_\Psi B^*-\Gamma(AWA^*+W^*)\right)(\la I_n-A^*)^{-1}C^*.
\]
Thus, to prove our claim (and because $(A^*,C^*)$ is reachable), we need to show that
\[
D_\Psi B^*-\Gamma(AWA^*+W^*)
\]
vanishes.
Substituting the value for $D_\Psi$ from (\ref{Dpsi}) in the above the expression we get
\begin{eqnarray*}
&&-\Gamma(H^*B^*-\bR B(B^*B)^{-1}B^*+AWA^*+W^*\\
&=&-\Gamma(W-\bR B(B^*B)^{-1}B^*+W^*)\\
&=&-\Gamma\bR(I_n-B(B^*B)^{-1}B^*).
\end{eqnarray*}
Recall that $\Gamma=\lim_{\epsilon\to 0}\Gamma_\epsilon$,
from the proof of Proposition \ref{lemma0}, while $\Gamma_\epsilon$ satisfies
\[
\Gamma_\epsilon (\bR+\epsilon \bPi_{\cN(\bR)}) =\Omega_\epsilon^{-1}B^*.
\]
Thus
\[
\Gamma_\epsilon (\bR+\epsilon \bPi_{\cN(\bR)})(I_n-B(B^*B)^{-1}B^*)=0
\]
identically for all $\epsilon$, and hence, taking the limit as $\epsilon\to 0$ we get the desired conclusion. This completes the proof of claim (i).

We first argue that $F_{ME}(\la)$ is analytic in $\mD$. Of course, $\Psi(\la)$ is already analytic in $\mD$ by our standing assumption on the location of the eigenvalues of $A$.
(Its poles cancel with the corresponding zeros of $\Phi(\la)^{-1}$ anyway.)
We only need to consider $\Phi(\la)^{-1}$. If $\bR$ is invertible, then $\Phi(\la)^{-1}$ has no poles in $\mD$ by \cite[Proposition 1]{acmatrix2}. If $\bR$ is singular, then, once again, we consider
\[
\bR_\epsilon=\bR+\epsilon \bPi_{\cN(\bR)}, \mbox{ with }\epsilon>0.
\]
With $\Omega_\epsilon=(B^*\bR_\epsilon^{-1}B)^{-1}$ and $\Gamma_\epsilon=\Omega_\epsilon B^*\bR_\epsilon^{-1}$ as before we define $\Phi_\epsilon(\la):=\Gamma_\epsilon G(\la)$ and apply \cite[Proposition 1]{acmatrix2} to deduce that $\Phi_\epsilon(\la)^{-1}$ is analytic in the closed unit disc, for all $\epsilon>0$. By continuity, $\Phi(\la)$ has no poles in the open unit disc. Similarly, the Hermitian part of $F_{ME}(\la)$ in $\mD$ is the limit of the Hermitian part of
\[
F_{ME,\epsilon}(\la):=\Phi_\epsilon(\la)^{-1}\Psi_\epsilon(\la)
\]
where $\Psi_\epsilon(\la)$ is given by (\ref{Psi}) with $\Gamma,\bR$ replaced by $\Gamma_\epsilon,\bR_\epsilon$, respectively. A matricial version of a classical identity between orthogonal polynomials (of first and second kind \cite[equation (1.17)]{Geronimus}) holds here as well:
\begin{equation}\label{PhiPsi}
\Psi(\la)\Phi(\la)^*+\Phi(\la)\Psi(\la)^*=\Omega.
\end{equation}
To verify this, after standard algebraic re-arrangement, the left hand side becomes
\[
\Lambda_0+\la \Gamma(I_n-\la A)^{-1}B_++\la^{-1} B_+^*(I_n-\la^{-1}A^*)^{-1}\Gamma^*
\]
where
\begin{eqnarray*}
B_+&=&A(BD_\Psi^*-\bR\Gamma^*+BH\Gamma^*), \mbox{ and}\\
\Lambda_0&=&D_\Psi+D_\Psi^*-\Gamma A\bR A^*\Gamma^*.
\end{eqnarray*}
If $\bR$ is invertible it is straightforward to show that
\[
BD_\Psi^*-R\Gamma^*+BH\Gamma^*=O_{n,m}
\]
while
\[
\Lambda_0=(B^*\bR^{-1}B)^{-1}=\Omega.
\]
If $\bR$ is singular then, as usual, we replace $\Gamma,\bR$ by their $\epsilon$-perturbations and claim the same identities for the relevant limits.
This shows that $F_{ME,\epsilon}(\la)\in\cF$ for all $\epsilon>0$. Hence, so is $F_{ME}(\la)$ since it is analytic in $\mD$ and its Hermitian part is nonnegative being the limit of the Hermitian part of $F_{ME,\epsilon}(\la)$ as $\epsilon\to 0$.
\end{proof}

\begin{remark} The relationship (\ref{PhiPsi}) (cf.\ \cite[equation (1.17)]{Geronimus})
between matricial functions of the ``first'' and ``second-kind'' generalizes to a two-sided version. Indeed, if we introduce analogous quantities for a right fraction
\[
F_{ME}(\la)=\Psi_r(\la)\Phi_r(\la)^{-1},
\]
by taking
$\Phi_r(\la)$
as in (\ref{Phir}) and
\begin{eqnarray*}
\Psi_r(\la)&:=&-L^*zA(I_n-zA)^{-1}\Gamma_r+D_{\Psi_r}\\
D_{\Psi_r}&=&=-(CC^*)^{-1}C(C^*L^*-\bR)\Gamma_r,
\end{eqnarray*}
then these satisfy
\begin{eqnarray}\nonumber 
&&\left[\begin{array}{cc} \Psi_\ell(\la)& \Phi_\ell(\la)\\
\Phi_r(\la)^* &\Psi_r(\la)^*
\end{array}\right]
\left[\begin{array}{cc} \Phi_r(\la) & \Phi_\ell(\la)^*\\ -\Psi_r(\la)&\Psi_\ell(\la)^* \end{array}\right] = \\
\label{PhiPsileftright} &&\hspace*{-25pt}
\left[\begin{array}{cc} O_m &\Omega_\ell\\ \Omega_r &O_m \end{array}\right]:=
\left[\begin{array}{cc} O_m &(C\bR^\sharp C^*)^\sharp\\ (B^*\bR^\sharp B)^\sharp &O_m \end{array}\right].
\end{eqnarray}
In the above we subscribe $\ell$, setting $\Phi_\ell(\la)=\Phi(\la)$ and $\Psi_\ell(\la)=\Psi(\la)$, to highlight ``left functions'' since $\Phi(\la),\Psi(\la)$ are the entries of the left fraction $F_{ME}(\la)=\Phi(\la)^{-1}\Psi(\la)$ of $F_{ME}(\la)$.
$\Box$
\end{remark}

\section{\bf Multivariable ``residues'' and singular parts}\label{residues}

We begin with
\begin{eqnarray}\label{Fme}
F(\la)&:=&F_{ME}(\la) =\Phi(\la)^{-1}\Psi(\la)\\
&=&D_\Psi+C_o\la(I_n-\la A_o)^{-1}B_o\nonumber
\end{eqnarray}
as given in (\ref{reducedFME}), suppressing the subscript ``$ME$'' for convenience.
When $\bR>0$, then $\Phi(\la)$ remains invertible in the closed unit disc and (\ref{PhiPsi}) readily implies that
\begin{equation}\label{matching}
\Herm\{F(e^{j\theta})\}=\Phi(e^{j\theta})^{-1}\Omega \left(\Phi(e^{j\theta})^{-1}\right)^*,
\end{equation}
cf.\ (\ref{mudot}). But when $\bR$ is singular, the variance of the minimal prediction error $\Omega$ is also singular (see Proposition \ref{Omegasingular} below) and (\ref{mudot}) may no longer be valid. The boundary limit of the Hermitian part defines a measure which may no longer be absolutely continuous. However, because $F(\la)$ is rational the singular part consists of finitely many disconinuities in $\mu(\theta)$. In order to separate the singular part from the absolutely continuous, we need to isolate the boundary poles of $F(\la)$. Accordingly, $F(\la)$ decomposes into a sum of ``lossless'' and ``lossy'' components---the lossless part being responsible for the singular part of the measure.

In the case where $F(\la)\in\cF$ is scalar-valued, the multiplicity of any pole
\[
\xi\in\partial \mD:=\{\la\;:\;\|\la\|=1\}
\]
cannot exceed one and $F(\la)$ decomposes into
\[
\rho\left(\frac{1+\la/\xi}{1-\la/\xi}\right) + F_{\rm remaining}(\la) \mbox{ with }\rho>0,
\]
where the first term is ``lossless'' and the second, $F_{\rm remaining}(\la)\in\cF$, has no singularity at $\xi$.
Conformably,
\[
d\mu(\theta)=\rho \,d\mU(\theta-\varangle\xi)+d\mu_{\rm remaining}(\theta)
\]
where $\varangle\xi$ denotes the angle of $\xi$ (i.e., $\xi=e^{j\varangle\xi}$) and $d\mu_{\rm remaining}(\theta)$ is continuous at $\varangle\xi$. Thus, in general,
\[
F(\la)=\sum_{\ell=1}^q \rho_i\left(\frac{1+\la/\xi_i}{1-\la/\xi_i}\right) + F_{\rm lossy}(\la)
\]
and the corresponding measure
\[
d\mu(\theta)=\sum_{i=1}^q \rho_i\,d\mU(\theta-\theta_i)+\dot{\mu}(\theta)d\theta.
\]

Analogous facts hold true in the multivariable case with some exceptions. Singularities in $\bR$ may not necessarily be associated with discontinuities in the measure and, while $F(\la)$ can have poles with higher multiplicity on the boundary of $\mD$, these may not have geometric multiplicity exceeding one. When $F(\la)$ has poles on the boundary, these are associated with discontinuities and our interest is to show how to decompose $F(\la)$ into a lossless and a lossy part, in general, and thus isolate the singular part of the measure. 

We first discuss the significance of $\bR$ being singular.
With $A_o,B_o,C_o$ as in (\ref{reducedMEdata}) and $\Gamma,\Omega$ as in Proposition \ref{lemma0} it holds that
\begin{equation}\label{RBOmegaB}
\bR=B\Omega B^*+A_o\bR A_o^*.
\end{equation}
This can be verified directly (by careful algebra).
It can also be shown via a limiting argument, replacing $\bR,\Gamma,A_o$ with $\bR_\epsilon,\Gamma_\epsilon,(I_n-B\Gamma_\epsilon)A$ (as in the proof of Theorem \ref{thm6})
and invoking \cite[equation (23)]{acmatrix2} to show that a similar identity holds for the perturbed quantities for all $\epsilon>0$, hence for their limits as well. A direct consequence of (\ref{RBOmegaB}) is the following.

\begin{prop}\label{Omegasingular}{\sf Let $A,B$ satisfy (\ref{standing1}a-d), $\bR\geq 0$, $A,B,\bR$ satisfy (\ref{rankcondition}),
and $\Omega$ the $\mH_m$-minimal value of $q_\bR$ subject to (\ref{GammaBI}). If $\Omega>0$ then $\bR>0$.
}\end{prop}
\begin{proof} The pair $(A_o,B\Omega^{1/2})$ is a reachable pair since it is obtained from $(A,B)$ after a state-feedback transformation and an invertible input tranformation. Then $\bR$ must be the reachability Grammian from (\ref{RBOmegaB}) which cannot be singular.
\end{proof}

\begin{ex}\label{cases}
Elementary scalar examples suffice to demonstrate how singularities of $\bR$ can give rise to poles of $F(\la)$ on $\partial \mD$. To see that this may not always be the case consider $A,B$ as in (\ref{companion}) with $n=4$ and $m=2$, and let $\bR$ which is now block-Toeplitz as in (\ref{bRsubn}) have entries
\[
R_0=\left[\begin{array}{cc}1&1\\1&1\end{array}\right] \mbox{ and }R_1=\frac{1}{2}R_0.
\]
Then
\[
\Gamma=\left[\begin{array}{cc}I_2&-\frac{1}{4}R_0\end{array}\right] \mbox{ and }\Omega=\frac{3}{4}R_0.
\]
Both $\bR$ and $\Omega$ are singular while the eigenvalues of $A_o$ are $\{0,0,0,\frac{1}{2}\}$. $\Box$
\end{ex}

Next, we present some general facts about lossless rational matrices in $\cF$. If
\[
d\mu(\theta)=\sum_{\ell=1}^q V_\ell \rho_\ell V_\ell^* d\mU(\theta-\theta_\ell)
\]
then
\[
\cH[d\mu(\theta)]=D_s+C_s\la(I_{n_s}-\la A_s)^{-1}B_s=:F_s(\la)
\]
where
\begin{eqnarray*}
D_s&=& \sum_{\ell=1}^q V_\ell \rho_\ell V_\ell^*\\
C_s&=& 2\left[\begin{array}{ccc}e^{j\theta_1}V_1&\ldots&e^{j\theta_q}V_q\end{array}\right]\\
B_s&=& \left[\begin{array}{c}\rho_1V_1^*\\\vdots\\\rho_qV_q^*\end{array}\right]\
\end{eqnarray*}
and $A_s$ block diagonal with blocks of the form $e^{j\theta_\ell} I_{n_\ell}$ of size equal to the size of $\rho_\ell$.
Then $F_s(s)\in\cF$ but it is also {\em lossless}, which amounts to $\Herm\{F_s(re^{j\theta})\}=0$ a.e.\ on $\partial \mD$.
It is a consequence of the Herglotz representation that,  modulo a state transformation and an additive skew-Hermitian summand in $D$, any rational lossless function is necessarily of this form.
An alternative characterization of lossless functions can be obtained via the well-known positive real lemma (e.g., \cite{FCG}) which, for the case where the Hermitian part is to be identically zero, specializes to the following.

\begin{prop}{\sf A rational function $D+C\la (I_n-\la A)^{-1}B$ belongs to $\cF$ and has Hermitian part identically equal to zero a.e.\ on the boundary of the unit circle if and only if there exists $P\geq 0$ such that
\begin{eqnarray}\label{conditionspra}
P-A^*PA&=&0\\\label{conditionsprb}
C^*-A^*PB&=&0\\\label{conditionsprc}
D+D^*-B^*PB&=&0.
\end{eqnarray}
}\end{prop}

\begin{proof}Nonnegativity of
\begin{equation}\label{matrixpr}
\left[ \begin{array}{cc}P-A^*PA &C^*-A^*PB\\ C-B^*PA&D+D^*-B^*PB\end{array}\right]
\end{equation}
along with $P\geq 0$ is equivalent to $D+C\la (I_n-\la A)^{-1}B\in\cF$ by the positive real lemma (see \cite[page 70]{FCG}).
Now consider its Hermitian part
\begin{eqnarray*}
\left[\begin{array}{cc}B^*\la^{-1}(I_n-\la^{-1}A^*)^{-1}&I_m\end{array}\right]\\
\times\left[\begin{array}{cc}O_n&C^*\\C&D+D^*\end{array}\right]
\left[\begin{array}{c}\la(I_n-\la^{-1}A)^{-1}B\\I_m\end{array}\right]
\end{eqnarray*}
and note that the null space of the mapping
\begin{eqnarray*}
M\mapsto {\mathcal{G}}(\la)^*M{\mathcal{G}}(\la),
\end{eqnarray*}
where
\[ {\mathcal{G}}(\la)=\left[\begin{array}{c}\la(I_n-\la A)^{-1}B\\I_m\end{array}\right],
\]
consists of matrices of the form
\[\left[ \begin{array}{cc}P-A^*PA &-A^*PB\\ -B^*PA&-B^*PB\end{array}\right].
\]
It readily follows that if conditions (\ref{conditionspra}-\ref{conditionsprc}) hold, then the function is lossless. If on the other hand (\ref{conditionspra}-\ref{conditionsprc}) do not hold and (\ref{matrixpr}) is simply nonnegative but not zero, then it can be shown that the Hermitian part can be factored into the product of nonzero spectral factors (cf.\ \cite[page 125]{FCG}). 
\end{proof}

Returning to (\ref{Fme}), in case $A_o$ has all its eigenvalues in the open disc $\mD$, then (\ref{mudot}) is valid and (\ref{matching}) holds as well for all $\theta$.
In case $A_o$ has eigenvalues on $\partial \mD$, we need to decompose $F(\la)$ into a lossless and a lossy summands.
To do this,
select $T_1,T_2$ matrices whose vectors form bases for the eignespaces of $A$ corresponding to eigenvalues on $\partial \mD$ and those in the interior of the disc, respectively. Then
$A_o$ transforms into a block triangular matrix
\begin{eqnarray*}
T^{-1}A_oT&=&\left[\begin{array}{cc}A_1 & 0\\0 & A_2\end{array}\right]
\end{eqnarray*}
where the spectrum of $A_1$ is on the boundary and of $A_2$ in the interior of the unit disc, respectively.
The input and output matrices $B_o,C_o$ transform conformably into
\begin{eqnarray*}
T^{-1}B_o&=&\left[\begin{array}{c}B_1 \\B_2\end{array}\right]\\
C_oT&=&\left[\begin{array}{cc}C_1 &C_2\end{array}\right].
\end{eqnarray*}
and
\[
F(\la)=D_\Psi + C_1\la(I -\la A_1)^{-1}B_1 +C_2\la(I -\la A_2)^{-1}B_2.
\]
Then we need to determine a value for a constant $D_1$ so that
\[F_1(\la)=D_1+C_1\la(I -\la A_1)^{-1}B_1
\]
is lossless. Necessarily, the remaining term $D_\Psi-D_1+C_2\la(I -\la A_2)^{-1}B_2$ is in $\cF$ and is devoid of singularities on the boundary.

The transformation $T_1$ above, can be chosen so that $A_1$ is unitary, since $A_1$ has only simple eigenvalues on $\partial \mD$.
Then condition (\ref{conditionspra}) leads to
\[
A_1P=PA_1
\]
and hence that $P$ is a polynomial function of $A_1$, i.e.,
\[P=p(A_1):=p_0 I +p_1 A_1+\ldots +p_{n_1-1} A_1^{n_1-1},
\]
$n_1$ being the size of $A_1$. The vector of coefficients $[\begin{array}{ccc}p_0&\ldots&p_{n_1-1}\end{array}]$ can now be computed from
(\ref{conditionsprb}) which becomes
\[A_1C^*=p(A_1)B_1.
\]
When $m>1$, this is an overdetermined set of equations which necessarily has a solution. Finally, we may take
\[
D_1=\frac{1}{2}B^*p(A_1)B
\]
to satisfy (\ref{conditionsprc}) and ensure that $F_1(\la)$ is lossless.
The matricial residues which represent the discontinuities in $d\mu(\theta)$ can now be computed by taking suitable limits at the singularities of $A_1$
\[
V_\ell \rho_\ell V_\ell^*=\Herm\{\lim_{\la\to e^{j\theta_\ell}} (1-\la e^{j\theta_\ell})F_1(\la)\},\;\ell=1,2,\ldots
\]
Evidently, if $A_1$ is first brought into a diagonal form, then a convenient closed expression for the limit can be given in terms of partitions of $B_1,C_1$ corresponding to the eigenvalue $e^{j\theta_\ell}$.

\section{\bf Impossibility of decomposition into white noise $+$ deterministic part}\label{whitepluscolor}

For the case of a scalar stochastic process $\{u_k\;:\; k\in\mZ\}$, where $m=1$, any state-covariance $\bR$ can be written as
\[\bR=\bR_{\rm signal} + \bR_{\rm white\; noise}
\]
where
\[
\bR_{\rm white\; noise} = \alpha_0 \bR_0
\]
with $\bR_0$ being the solution to the Lyapunov equation
\[
\bR_0-A\bR_0 A^*=BB^*
\]
and $\alpha_0$ the
smallest eigenvalue of the matrix pencil $\bR-\alpha \bR_0$, i.e.,
\begin{eqnarray}\label{interpretationa}
\alpha_0&=&\min\{\alpha\;:\; \det\left(\bR-\alpha \bR_0\right)=0\}\\\label{interpretationb}
&=&\max\{\alpha\;:\;\bR-\alpha \bR_0\geq 0\}.
\end{eqnarray}
The matrix $\bR_0$ is the controllability Grammian of the pair $(A,B)$ and represents the state-covariance when the input is unit-variance white noise. Then $\bR_{\rm white\; noise}$ represents the maximal summand of $\bR$ that can be attributed to a white-noise input component of (\ref{i2s}), while the remaining $\bR_{\rm signal}$ corresponds to a deterministic input part. It can also be shown that this decomposition is canonical in the sense that any other one, consistent with a ``white noise plus deterministic part'' hypothesis for the input, will have a larger number of deterministic components (i.e., spectral lines).
This is the interpretation of the 
CFP decomposition.
The theory was originally developed for $\bR$'s having a Toeplitz structure \cite{GrenanderSzego,StoicaMoses} and extended to general state-covariances in \cite{sesha1,sesha2}.

It is rather instructive to present a derivation of the fact that, when $m=1$, the equivalent
conditions (iii, iii-a) of Proposition \ref{lemma0} are automatically satisfied by any singular state-covariance. This underscores the dichotomy with the multivariable case where a decomposition of $\bR$ consistent with a ``white noise plus deterministic part'' input is not always possible (see Examples \ref{example1} and \ref{example2} below).

\begin{prop}\label{mequalone}{\sf Let $\bR,A,B,H$ satisfy (\ref{rankcondition}b) in Theorem \ref{thm2}, let $\bR\geq 0$ and singular, and let $m=1$. Then $B^*\bPi_{\cN(\bR)}B$ is invertible. 
}\end{prop}
\vspace*{.1in}

\begin{proof} Suppose that $B^*\bPi_{\cN(\bR)}B$ is not invertible. Then
\begin{equation}\label{RBequal0}
\bPi_{\cN(\bR)}B=O_{n\times 1}
\end{equation}
and
\begin{equation}\label{containment}
\cR(B)\subseteq \cR(\bR).
\end{equation}
From (\ref{rankcondition}b) and (\ref{RBequal0}) it follows that
$\bPi_{\cN(\bR)}A\bR A^*\bPi_{\cN(\bR)} = O_{n\times n}$, and hence, that
\[\bPi_{\cN(\bR)}A\bR = O_{n\times n}.
\]
From (\ref{containment}), $\bPi_{\cN(\bR)}AB=O_{n\times 1}$.
By induction, using (\ref{rankcondition}b), it follows that
\[
\bPi_{\cN(\bR)}A^\ell \bR=O_{n\times n}, \mbox{ for }\ell=0,1,\ldots
\]
and hence, that $\cR(\bR)$ is $A$-invariant. But $\cR(B)\subseteq \cR(\bR)$ and so is the largest $A$-invariant subspace containing $\cR(B)$. Because $(A,B)$ is a reachable pair, $\cR(\bR)=\mC^n$ which contradicts the hypothesis that $\bR$ is singular.
\end{proof}

The following example shows that the statement of the proposition is only valid when $m=1$ and that, in general, a decomposition of $\bR$ consistent with a ``white noise plus deterministic part'' input is not always possible.

\begin{ex}\label{example1} Let
\[A=\left[ \begin{array}{cc}O_2& O_2\\ I_2 &O_2\end{array}\right],\; B=\left[ \begin{array}{c}
                                       I_2\\ O_2\end{array}\right]
\]
and
\[
\bR=\left[
\begin{array}{cccc} 1&0&1/2 &3/4\\
                                0&1&0 & 1/2\\
                      1/2 &0&1&0\\
                      3/4 & 1/2 &0 &1\end{array}
\right],
\]
where, as usual, $I_2$ and $O_2$ are the $2\times 2$ identity and zero matrices, respectively.
It can be readily seen that they satisfy conditions (\ref{standing1}) as well as (\ref{rankcondition}b) in Theorem \ref{thm2}---$\bR$ being a block-Toeplitz matrix.
Then $\bR\geq 0$ and singular. To see this note that the first three principal minors of $\bR$ are positive definite while
\[
\left(\begin{array}{cccc}-2&-1&1&2\end{array}\right)\bR= O_{1\times 4}.
\]

If the input to (\ref{i2s}) is white noise with variance the $2\times 2$ non-negative matrix
\[
Q=\left[\begin{array}{cc} a& b\\ \bar{b}& c\end{array}\right],
\]
then the state-covariance (for the chosen values of $(A,B)$ and corresponding to this white-noise input) is
\[
\bR_0=\left[\begin{array}{cc} Q& O_2\\ O_2& Q\end{array}
\right] = I_2\otimes Q.
\]
We claim that
\[
\bR-\bR_0\geq 0 \Rightarrow \bR_0 =O_{4\times 4}.
\]
To prove this, consider that $Q\geq 0$ from which we obtain
\begin{equation}\label{Qgeq0}
ac\geq |b|^2,\; a\geq 0,\; c\geq 0.
\end{equation}
Now, if $v=\left(\begin{array}{cccc}-2&-1&1&2\end{array}\right)$ then $v\bR v'=0$ and $v \bR_0 v'\geq 0$.
Therefore
\begin{eqnarray}\nonumber
&&\bR-\bR_0 \geq 0\\\nonumber
&\Rightarrow& v\bR_0 v'= 0\\\nonumber
&\Rightarrow& \left(\begin{array}{cc}2&1\end{array}\right)Q \left(\begin{array}{c}2\\1\end{array}\right)+
\left(\begin{array}{cc}1&2\end{array}\right)Q \left(\begin{array}{c}1\\2\end{array}\right)= 0\\
&\Rightarrow& 5a+4 \Real(b)+5c=0.\label{algebraic}
\end{eqnarray}
Thence, if $\beta:=\Real(b)$,
\begin{eqnarray*}
&&ac\geq \beta^2\\
&\Rightarrow& ac\geq \frac{25}{16}(a^2+c^2+2ac)\\
&\Rightarrow& 0\geq a^2 + c^2 + \frac{34}{25}ac\\
&\Rightarrow& \mbox{ either } a=0 \mbox{ or } c=0.
\end{eqnarray*}
In either case, $|b|=0$ and hence all three $a=b=c=0$ from (\ref{algebraic}).
Thus, $Q=O_{2\times 2}$ and $\bR_0 =O_{4\times 4}$ as claimed. $\Box$
\end{ex}

While the previous example shows that no white noise component can be subtracted
in the hope of reaching a state-covariance satisfying condition (iii) in Proposition \ref{lemma0} (thus corresonding to pure sinusoids), more is true. The following example shows that the off-diagonal block-entries of a block-Toeplitz $\bR$ already prevent condition (iii) from being true.

\begin{ex}\label{example2}Let $A,B$ as in Example \ref{example1} and
\[
\bR=\left[
\begin{array}{cccc} a&b&1/2 &3/4\\
                               \bar{b}&c&0 & 0\\
                      1/2 &0&a&b\\
                      3/4 & 0 & \bar{b}&c\end{array}
\right].
\]
In order for condition (iii) of Proposition \ref{lemma0} to hold, the null space $\cN(\bR)$ must have a dimension $\geq 2 = {\rm dim}(\cR(B))$ (which can also readily seen from condition (iii-a) as well). We argue that this cannot happen.
Since
\[
\left[\begin{array}{cc} a&3/4\\
                               3/4&c\end{array}\right]
\]
is a principle minor of $\bR\geq 0$, neither $a$ nor $c$ can vanish.
The rank of
\[
R_0:=\left[\begin{array}{cc} a&b\\
                               \bar{b}&c\end{array}
\right]
\]
must be equal to one, since there is a $3\times 3$ minor of $\bR$ with determinant
\[
c\times \det (R_0).
\]
Hence, $ac=|b|^2 \Rightarrow c=|b|^2/2a$. 
But then, the northwest $3\times 3$ principle minor of $\bR$ is equal to
$-c/4<0$, which contradicts $\bR\geq 0$.
$\Box$
\end{ex}

\section{\bf Decomposition as a convex optimization problem}\label{shortrange}

We have just seen that in the case of a vectorial input, a decomposition of the state-covariance $\bR$ of (\ref{i2s}) which is consistent with the hypothesis of ``white noise plus a deterministic signal at the input'' may not always be possible. We begin by choosing an alternative interpretation of the
CFP decomposition as seeking to separate the maximal-variance white noise component at the input which is consistent with a known state-covariance. This is the analog of (\ref{interpretationb}) and leads to the following problem.

\begin{problem}\label{mainthm3}{\sf Given $\bR,A,B$ satisfying (\ref{standing1}), $\bR\geq 0$, and (\ref{rankcondition}a) in Theorem \ref{thm2},
determine a decomposition
\begin{equation}\label{decomposition_problem1}
\bR=\bR_{\rm signal}+\bR_{\rm noise}
\end{equation}
where the summands satisfy
\begin{eqnarray}\label{nonegative1}
\bR_{\rm noise}&\geq & 0,\\\label{nonegative2}
\bR_{\rm signal}&\geq & 0,\\\label{condition3}
\hspace*{-25pt}\bR_{\rm noise}-A\bR_{\rm noise} A^*&=&B Q B^* \mbox{ with }Q\geq 0,
\end{eqnarray}
and
\begin{equation}\label{argmax}
\bR_{\rm noise}={\rm argmax}\{ {\rm trace\,}\bR_{\rm noise}\;:\; (\ref{nonegative1}-\ref{condition3}) \mbox{ hold}\}.
\end{equation}
}\end{problem}
\vspace*{.1in}

This is a standard convex optimization problem where the noise variance ${\rm trace\,}Q$ is a linear functional of the parameters in $Q$ and all constraints appear in the form of linear matrix inequalities. Thus, it can be readily and efficiently solved with existing computational tools.
Alternatives to (\ref{argmax}) corresponding to a different ``normalizations'' are
\begin{equation}\label{argmax2}
\bR_{\rm noise}={\rm argmax}\{ {\rm trace\,}\left(\bR_{\rm noise}{\bf W}\right)\;:\; (\ref{nonegative1}-\ref{condition3}) \mbox{ hold}\},
\end{equation}
for any weight matrix ${\bf W}>0$ (which may encapsulate ``prior'' information about the directionality of the noise), or to seek
\begin{equation}\label{argmax3}
Q={\rm argmax}\{ {\rm trace\,}Q\;:\; (\ref{nonegative1}-\ref{condition3}) \mbox{ hold}\}.
\end{equation}

Below we present an example which shows that a maximum-trace solution as above, in general, does not lead to a decomposition with $\bR_{\rm signal}$ corresponding to a deterministic signal (i.e., satisfying (\ref{singleton})) even when an alternative decomposition does.

\begin{ex} With $A,B$ as in Example \ref{example1}, consider the state-covariance
\[
\bR=\left[
\begin{array}{cccc} r_0&r_1&1/2 &3/4\\
                               r_1&r_0&0 & 1/2\\
                      1/2 &0&r_0&r_1\\
                      3/4 & 1/2 &r_1&r_0\end{array}
\right],
\]
where the block-diagonal entries are yet unspecified. The values for these entries can be explicitly computed in the following two cases:
\begin{itemize}
\item[(i)] $B^*\bPi_{\cN(\bR)}B$ is invertible, and
\item[(ii)] ${\rm trace}(\bR)$ is minimal,
\end{itemize}
while always $\bR\geq 0$.

The first can be carried out as follows. Condition (i) is equivalent to the existence of
a matrix
\[
\Gamma=\left[\begin{array}{cccc} 1&0& \gamma_{1,3} &\gamma_{1,4}\\
                                                      0&1&\gamma_{2,3}&\gamma_{2,4}
                                                      \end{array}\right]
\]
such that $\Gamma \bR$ is the zero matrix. Denote
\begin{eqnarray*}
R_0&:=&\left[\begin{array}{cc} r_0&r_1\\r_1&r_0\end{array}\right],\; 
                                                  R_1:=    \left[
\begin{array}{cc} 1/2 &3/4\\
                            0 & 1/2\end{array}
\right], \mbox{ and}\\
\Gamma_0&:=&\left[\begin{array}{cc} \gamma_{1,3} &\gamma_{1,4}\\
                                                      \gamma_{2,3}&\gamma_{2,4}
                                                      \end{array}\right].
\end{eqnarray*}
Since
\begin{eqnarray*}
R_0+\Gamma_0 R_1^*=O_2 &\Rightarrow & R_0 +R_1\Gamma_0^*=O_2, \mbox{ while}\\
R_1+\Gamma_0 R_0=O_2,
\end{eqnarray*}
we deduce that
\begin{eqnarray}
R_1-\Gamma_0 R_1 \Gamma_0^*&=&O_2, \label{first}\\
R_1-\Gamma_0^2 R_1^*&=&O_2.\label{second}
\end{eqnarray}
Equation (\ref{first}) leads to
\[
R_1+R_1^*=\Gamma_0 (R_1+R_1^*) \Gamma_0^*
\]
and, if we factor $R_1+R_1^*=SS^*$ with
\[S=\left[
\begin{array}{cc} 1 &0\\
                            3/4 & \sqrt{1-\left(\frac{3}{4}\right)^2}\end{array}
\right],
\]
we deduce that $S^{-1}\Gamma_0S$ must be unitary. Then from (\ref{second})
we determine the eigenvalues of $\Gamma_0$. Carrying out all computations explicitely leads to
\[
R_0=\frac{1}{2}\left[
\begin{array}{cc} \frac{1}{\cos(\theta)} &\tan(\theta)\\
                            \tan(\theta) &  \frac{1}{\cos(\theta)} \end{array}
\right]
\]
and
\[
\Gamma_0=\left[
\begin{array}{cc} \frac{\cos(2\theta)}{\cos(\theta)} &-\tan(\theta)\\
                            \tan(\theta) &  \frac{1}{\cos(\theta)} \end{array}
\right]
\]
where $\sin(\theta)=\frac{3}{4}$. The values in $R_0$ is the unique set values for which (i) holds.

Similarly, the computation of the state-covariance with minimal trace as in (ii) can be carried out explicitly to give
\[R_{0,\rm min\,trace}=\left[
\begin{array}{cc} 3/4 &1/2\\
                        1/2&    3/4\end{array}
\right].
\]
Finally, it is easy to check that $R_0- R_{0,\rm min\,trace}$ is indefinite. $\Box$
\end{ex}



\section{\bf Short-range correlation structure}\label{sec:shortrange}

The rationale for the 
CFP decomposition has been re-cast in Problem \ref{mainthm3} as seeking to extract the maximal variance that can be attributed to white-noise. In the case where $\bR$ is block-Toeplitz as in (\ref{bRsubn}), this amounts to determining a block-diagonal matrix $\bR_{\rm noise}$ of maximal trace satisfying the required positivity constraints (\ref{nonegative1}-\ref{condition3}). Yet, it is rarely the case in practice that a ``white-noise'' hypothesis is valid. Thus, we herein propose a new paradigm--a paradigm that also leads to a convex optimization problem and encompasses the above interpretation of the 
CFP decomposition as a special case. We seek to identify a maximal-variance summand which has a ``short-range correlation structure'' defined as follows:

\begin{definition}Given $A,B$ satisfying (\ref{standing1}) a state-covariance $\bR$ of the system (\ref{i2s}) has {\em correlation range} $k$ if
there exists a matrix $H\in\mC^{m\times n}$  so that
\begin{equation}\label{soughtH}
H^*=[\begin{array}{cccc}B&AB&\ldots &A^kB\end{array}]\left[\begin{array}{c}Q_0^*\\Q_1^*\\\vdots\\Q_k^*\end{array}\right]
\end{equation}
for suitable matrices $Q_0,\ldots,Q_k$, such that
\begin{equation}\label{soughtHrhs}
\bR-A\bR A^*=BH+H^*B^*
\end{equation}
and
\begin{equation}\label{needed}
Q_0+\la Q_1+\ldots+\la^kQ_k
\in\cF.
\end{equation}
\end{definition}
\vspace*{.1in}

It is insightful to first consider the case where 
$A,B$ are given as in (\ref{companion}) and the state-covariance structure is $(\ell+1)\times(\ell+1)$ block-Toeplitz. A block-Toeplitz matrix $\bR$
has correlation range $k$ if it is block-banded with all entries beyond the $k$th one being zero and, most importantly, it remains a covariance matrix when extended with zero elements beyond the $\ell$th entry as well. This is equivalent to
$R_{k+1}=R_{k+2}=\ldots=R_\ell=\ldots=O_m$ being an admissible extension since already
\[
R_0+2\la R_1+\ldots +2\la^k R_k \in\cF
\] 
from (\ref{needed}) because $Q_i=R_i$ ($i=1,\ldots,k$) and $R_0=Q_0^*+Q_0$.

\begin{ex} The following elementary example helps illustrate the concept of bounded correlation range. Consider
the Toeplitz matrix
\[
\bR=\left[ \begin{array}{ccc} 1 & 1/2 & 1/3\\
                                       1/2& 1 & 1/2\\
                                       1/3& 1/2 &1\end{array}\right].
\]
We seek a Toeplitz noise-covariance summand of maximal trace with correlation range $1$, i.e., we seek
\[
\bR_{\rm noise}=\left[ \begin{array}{ccc} q_0 & q_1 & 0\\
                                       q_1& q_0 & q_1\\
                                       0& q_1 &q_0\end{array}\right]
\]
so that $\bR-\bR_{\rm noise}\geq 0$, and $q_0+2\la q_1\in\cF$.
Since $q_0+2\la q_1$ is only of degree one,  $q_0+2\la q_1\in\cF$ if and only if $|q_1|\leq q_0$. The solution turns out to be $q_0=2/3$ and $q_1=0.3097$. 

Instead, if we sought $\bR_{\rm noise}$ diagonal corresponding to white noise, the answer would have been $\bR_{\rm noise}=\min\{ {\rm eig}(\bR)\} \times I_3$. It can be easily checked that $\min\{{\rm eig}(\bR)\}=0.4402<q_0$. Thus, colored MA-noise allows a larger amount of energy to be accounted for. $\Box$
\end{ex}

Problem \ref{mainthm3} with condition (\ref{condition3}) replaced by
\begin{equation}\label{condition3a}
\bR _{\rm noise} \mbox{ having correlation range }k
\end{equation}
is also a convex optimization problem. In general, the positive-real constraint (\ref{needed}) can be expressed as a convex condition via the well-known positive-real lemma (e.g., see \cite{FCG}), and the maximizer of the trace can be
readily obtained with existing numerical tools (e.g., the Matlab LMI toolbox).

In the case (\ref{i2s}) has nontrivial dynamics, the right hand side of (\ref{soughtHrhs}) becomes
\begin{eqnarray}\label{Qs}\nonumber
BH+H^*B^*&=&A^kBQ_k^*B^*+\ldots +ABQ_1^*B^*\\
\nonumber &&  \hspace*{-3cm}+ B(Q_0^*+Q_0)B^*+BQ_1B^*A^*+\ldots+BQ_kB^*(A^*)^k,
\end{eqnarray}
and $\bR$ can be interpreted as the state covariance due to colored noise at the input with spectral density
\[ Q_k^*e^{-j\theta}+\ldots +Q_1^*e^{-j\theta}+(Q_0^*+Q_0) + Q_1e^{j\theta} + \ldots +Q_k e^{jk\theta}.
\]
A detailed study on the potential of decomposition according to ``correlation range'' for high resolution spectral analysis will be presented in a forthcoming report.

\section{\bf Concluding remarks}
The Carath\'{e}odory-Fej\'{e}r-Pisarenko (CFP) decomposition underlies many subspace identification techniques in modern spectral analysis (such as MUSIC, ESPRIT, and their variants \cite{StoicaMoses}). But in spite of its importance and its extensive appearance in many guises in the identification and signal processing literature, no multivariable analog had been proposed. Perhaps the reason can be sought in the fact that the exact analog of the CFP-decomposition does not exist. This realization led us to alternative interpretations of the CFP-decomposition, and the goal of this paper has been to explore such alternatives for a ``signal plus noise'' decomposition of covariances for multivariable processes. 
In the process we have found that (e.g., see Example \ref{example2} and 
Section \ref{whitepluscolor}) regardless of how much of the energy is accounted for by noise, the remaining energy, in general, cannot be accounted for by pure spectral lines only. The remaining energy necessarily corresponds to a singular covariance matrix and thus, Sections \ref{positiverealfunction} and \ref{residues} develop the needed theory to construct spectra for singular matrices. Finally Sections \ref{shortrange} and \ref{sec:shortrange} develop certain alternatives to the CFP decomposition where we forgo the requirement that one part is completely deterministic, and allow instead that it has a long range correlation structure.




\section{Acknowledgments}

The author wishes to thank Dr.\ Dan Herrick for his input and for discussions that partially motivated this work.

\end{document}